# SOME LOCAL APPROXIMATIONS OF DAWSON–WATANABE SUPERPROCESSES

By Olav Kallenberg

*Auburn University*

Let $\xi$ be a Dawson–Watanabe superprocess in $\mathbb{R}^d$ such that $\xi_t$ is a.s. locally finite for every $t \geq 0$. Then for $d \geq 2$ and fixed $t > 0$, the singular random measure $\xi_t$ can be a.s. approximated by suitably normalized restrictions of Lebesgue measure to the $\varepsilon$-neighborhoods of supp $\xi_t$. When $d \geq 3$, the local distributions of $\xi_t$ near a hitting point can be approximated in total variation by those of a stationary and self-similar pseudo-random measure $\tilde{\xi}$. By contrast, the corresponding distributions for $d = 2$ are locally invariant. Further results include improvements of some classical extinction criteria and some limiting properties of hitting probabilities. Our main proofs are based on a detailed analysis of the historical structure of $\xi$.

**1. Introduction.** By a *Dawson–Watanabe superprocess* (or *DW-process*, for short) we mean a vaguely continuous, measure-valued Markov process $\xi$ on $\mathbb{R}^d$ satisfying $E_\mu \exp(-\xi_t f) = \exp(-\mu v_t)$ for any $f \in C_K^+(\mathbb{R}^d)$, where $v$ is the unique solution on $\mathbb{R}_+ \times \mathbb{R}^d$ to the *evolution equation* $\dot{v} = \frac{1}{2}\Delta v - v^2$ with initial condition $v_0 = f$. The more general process with $v^2$ replaced by $\frac{1}{2}\gamma v^2$ can be reduced to the present version by a suitable scaling. The usual construction for bounded initial measures $\mu$ extends by independence to any $\sigma$-finite initial measure $\mu$. By Lemma 3.2 below, $\xi_t$ is then a.s. locally finite for every $t > 0$ iff $\mu p_t < \infty$ for all $t$, where $p_t$ denotes the standard normal density $p_t(x) = (2\pi t)^{-d/2} \exp(-|x|^2/2t)$ on $\mathbb{R}^d$.

The DW-process has been studied intensely, along with more general superprocesses, for the last 30 years, and the literature on the subject is absolutely staggering with respect to both volume and depth. Several excellent surveys exist, including the lecture notes and monographs [3, 7, 8, 22, 26].

For $d \geq 2$ and a fixed $t > 0$, $\xi_t$ is known to be a.s. singular and diffuse with a support of Hausdorff dimension 2 (cf. Theorem 6.15 in [8]). Writing $\xi_t^\varepsilon$ for the restriction of Lebesgue measure $\lambda^d$ to the $\varepsilon$-neighborhood of $\operatorname{supp} \xi_t$ it was shown by Tribe [27] that $\varepsilon^{2-d} \xi_t^\varepsilon \xrightarrow{v} c_d \xi_t$ a.s. as $\varepsilon \to 0$ when $d \geq 3$, where $\xrightarrow{v}$ denotes vague convergence and $c_d > 0$ is a universal constant. For $d = 2$ we prove in Theorem 7.1 that $\tilde{m}(\varepsilon)|\log \varepsilon| \xi_t^\varepsilon \xrightarrow{v} \xi_t$ a.s., where the function $\tilde{m}$ is such that $\log \tilde{m}$ is bounded with strong continuity properties. In particular, this confirms that $\xi_t$ "distributes its mass over $\operatorname{supp} \xi_t$ in a deterministic manner" (cf. [8], page 115, or [26], page 212), as previously inferred from some deep results involving the exact Hausdorff measure (cf. [5]).

Our proofs depend crucially on some basic hitting estimates, due to Dawson, Iscoe and Perkins [4] for $d \geq 3$ and Le Gall [21] for $d = 2$. The former paper gives $\varepsilon^{2-d} P_\mu \{\xi_t B_0^\varepsilon > 0\} \to c_d \mu p_t$ for $d \geq 3$ as $\varepsilon \to 0$, where $B_x^r$ denotes an open ball around $x$ of radius $r$. Likewise, combining Le Gall's results with an analysis of the historical structure, we show in Theorem 5.3 that $\tilde{m}(\varepsilon)|\log \varepsilon| P_\mu \{\xi_t B_0^\varepsilon > 0\} \to \mu p_t$ for $d = 2$, with $\tilde{m}$ as before. A simple rescaling argument in Theorem 4.5 shows that the local extinction property $\xi_t \xrightarrow{d} 0$ as $t \to \infty$, first noted by Dawson [2] when $d = 2$ and $\xi_0 = \lambda^2$, is equivalent to the seemingly stronger support property $\operatorname{supp} \xi_t \xrightarrow{d} 0$. (Note that the two properties are given by $\xi_t B \xrightarrow{P} 0$ and $1\{\xi_t B > 0\} \xrightarrow{P} 0$, respectively, for any bounded Borel set $B$.)

Another main result is Theorem 8.1, where we show for $d \geq 3$ that the conditional distribution of $\xi_t$, given that $\xi_t$ charges a small set $B$, can be approximated in total variation by the corresponding conditional distribution for a certain stationary and self-similar pseudo-random measure $\tilde{\xi}$. (The prefix "*pseudo*" indicates that the underlying "probability" measure $\tilde{P}$ is not normalized and may even be unbounded. This anomaly is prompted by the self-similarity of $\tilde{\xi}$, as explained in [28]. In our context it causes no problems, since the associated hitting probabilities remain finite.) By contrast, we prove in Theorem 9.1 that for $d = 2$, the random measure $\xi_t$ is asymptotically invariant near a hitting point.

The present work is part of a general program outlined in [16], where we indicate how a whole class of local properties seem to be shared by three totally different types of random objects—by simple point processes, local time random measures, and certain measure-valued diffusion processes. The point process case is classical and has been thoroughly explored in [11, 17]. Versions of the Lebesgue approximation in Theorem 7.1 are known for the local time random measures of regenerative sets and exchangeable interval partitions (cf. [18] and Proposition 6.13 in [15]), and some delicate approximations related to Theorem 8.1 below appear in [12, 14].

As a referee points out, certain *intersection local time* random measures may be added to our list of random objects with related properties. For

APPROXIMATION OF SUPERPROCESSES 3

example, a Lebesgue approximation analogous to ours was proved in this context by Le Gall [20] (though with convergence in $L^p$, $p \geq 1$, rather than a.s.), and similar results on the *average density* have been obtained for DW-processes and intersection local times by Mörters and Shieh [23, 24], giving further evidence of the profound dichotomy between the cases when $d \geq 3$ or $d = 2$. It is also interesting to note that their results for the intersection local time are stated in terms of the Palm distribution associated with a suitable stationary and self-similar pseudo-random measure (cf. [24], pages 3f), corresponding to our $\tilde{P}^0$ in Theorem 8.2.

We proceed with some general remarks on terminology and notation. A *random measure* on $\mathbb{R}^d$ is defined as a locally finite kernel $\xi$ from some basic probability space $(\Omega, \mathcal{A}, P)$ into $(\mathbb{R}^d, \mathcal{B}^d)$, where $\mathcal{B}^d$ denotes the Borel $\sigma$-field on $\mathbb{R}^d$. Thus, $\xi(\omega, B)$ is a locally finite measure in $B \in \mathcal{B}^d$ for fixed $\omega \in \Omega$ and is measurable in $\omega$ for fixed $B$. A *pseudo-random* measure is defined in the same way, except that the underlying measure $\tilde{P}$ is now allowed to be $\sigma$-finite. We may also regard $\xi$ as a measurable function from $\Omega$ to the space $\mathcal{M}_d$ of locally finite measures on $\mathbb{R}^d$, equipped by the $\sigma$-field generated by all evaluation maps $\pi_B : \mu \mapsto \mu B$ with $B \in \mathcal{B}^d$. The subclasses of bounded sets and measures are denoted by $\hat{\mathcal{B}}^d$ and $\hat{\mathcal{M}}_d$, respectively.

The *vague* topology in $\mathcal{M}_d$ is generated by all integration maps $\pi_f : \mu \mapsto \mu f = \int f \, d\mu$ with $f$ belonging to the space $C_K^d$ of continuous functions $\mathbb{R}^d \to \mathbb{R}_+$ with bounded support. Similarly, the *weak* topology in $\hat{\mathcal{M}}_d$ is generated by the maps $\pi_f$ for all $f$ in the class $C_b^d$ of bounded, continuous functions $\mathbb{R}^d \to \mathbb{R}_+$. Thus, $\mu_n \xrightarrow{v} \mu$ in $\mathcal{M}_d$ iff $\mu_n f \to \mu f$ for all $f \in C_K^d$, and similarly for $\mu_n \xrightarrow{w} \mu$ in $\hat{\mathcal{M}}_d$. For random measures $\xi_n$ and $\xi$ on $\mathbb{R}^d$, the associated $L^1$-convergence $\xi_n \to \xi$ means that $\xi_n f \to \xi f$ in $L^1$ for all $f$ in $C_K^d$ or $C_b^d$, respectively. Convergence in distribution of random measures, denoted by $\xi_n \xrightarrow{d} \xi$, is understood to be with respect to the vague topology, unless something else is said. Note that this is equivalent to $\xi_n f \xrightarrow{d} \xi f$ for all $f \in C_K^d$ (cf. Theorem 16.16 in [13]).

Convergence of closed random sets is defined as usual with respect to the Fell topology (cf. [13], pages 324, 566). However, in this paper we need only the special cases of convergence to the empty set or the whole space, which are explained whenever they occur.

Throughout the paper we use relations such as $\underset{\simeq}{=}$, $\underset{\simeq}{\leq}$, $\underset{\simeq}{\simeq}$ and $\asymp$, where the first three mean equality, inequality and asymptotic equality up to a constant factor, and the last one is the combination of $\underset{\simeq}{\leq}$ and $\underset{\simeq}{\geq}$. We often write $a \ll b$ to mean $a/b \to 0$. The double bars $\|\cdot\|$ denote the supremum norm when applied to functions and total variation when applied to signed measures. We also write $\|\cdot\|_B$ for the supremum or total variation over the set $B$. For functions $f_n$ or signed measures $\mu_n$ on $\mathbb{R}^d$, the convergence $\|f_n\| \to 0$ or $\|\mu_n\| \to 0$ is said to hold *locally* if $\|f_n\|_B \to 0$ or $\|\mu_n\|_B \to 0$, respectively, for all $B \in \hat{\mathcal{B}}^d$. In Section 8 we also use the notation $\|\mu_n\|_B \to 0$ for signed measures $\mu_n$ on $\mathcal{M}^d$ and sets $B \in \hat{\mathcal{B}}^d$, in which case the precise meaning is explained in connection with Theorem 8.1.

In any Euclidean space $\mathbb{R}^d$, we write $B_x^r$ for the open ball of radius $r > 0$ centered at $x \in \mathbb{R}^d$. The shift and scaling operators $\theta_x$ and $S_r$ are given by $\theta_x y = x + y$ and $S_r x = rx$, respectively, and for measures $\mu$ on $\mathbb{R}^d$ we define $\mu\theta_x$ and $\mu S_r$ by $(\mu\theta_x)B = \mu(\theta_x B)$ and $(\mu S_r)B = \mu(S_r B)$, respectively. In particular, $(\mu S_r)f = \mu(f \circ S_r^{-1})$ for measurable functions $f$ on $\mathbb{R}^d$. Convolutions of measures $\mu$ with functions $f$ are given by $(\mu * f)(x) = \int f(x - u)\mu(du)$. Product measures are written as $\mu \otimes \nu$ or $\mu^n = \mu \otimes \cdots \otimes \mu$, and in particular $\lambda^d$ denotes Lebesgue measure on $\mathbb{R}^d$. The functional notations $f(x)$ and $f_x$ are used interchangeably, depending on typographical convenience. Notation pertaining to Palm measures or DW-processes is explained in the next section.

The paper is organized as follows. In Section 2 we prove some preliminary technical results and explain the crucial ideas about DW-processes, cluster representations and Palm measures needed in subsequent sections. In Section 3 we characterize the locally finite DW-processes in terms of their initial measures and derive some useful estimates of the second moments. In Section 4 we use the classical hitting estimates to give bounds on the associated multiplicities, and we establish some weak extinction criteria for $d \geq 2$. In Section 5 we identify and study the proper normalization for the hitting probabilities to converge when $d = 2$. In Section 6 we estimate the second moments of the neighborhood measures $\eta_t^\varepsilon$ associated with the clusters $\eta_t$ of a DW-process. In Section 7 we are ready to prove the mentioned Lebesgue approximation for DW-processes of dimensions $d \geq 2$. In Section 8 we prove the mentioned approximation in total variation for DW-processes of dimension $d \geq 3$. Finally, we show in Section 9 that DW-processes of dimension 2 are locally invariant in a number of different ways.

**2. Preliminaries.** In this paper DW-processes are often denoted by $\xi = (\xi_t)$, and we write $P_\mu\{\xi \in \cdot\}$ for the distribution of the process $\xi$ with initial measure $\mu$. The same notation is used for the entire historical process. In all the mentioned literature, $\xi$ is first constructed for bounded $\mu$. To extend the definition to the $\sigma$-finite case, we may write $\mu = \sum_n \mu_n$ for some bounded measures $\mu_n$, and choose $\xi_1, \xi_2, \ldots$ to be independent DW-processes starting from $\mu_1, \mu_2, \ldots$. Then $\xi = \sum_n \xi_n$ is a locally finite DW-process with initial measure $\mu$, provided that $\mu p_t < \infty$ for all $t > 0$.

For every fixed $\mu$, the DW-process $\xi$ is infinitely divisible under $P_\mu$ and admits a decomposition into a Poisson "forest" of conditionally independent *clusters*, corresponding to the excursions of the contour process in the ingenious "Brownian snake" representation of Le Gall [22]. In particular, this yields a cluster representation of $\xi_t$ for every fixed $t > 0$. More generally, the "ancestors" of $\xi_t$ at an earlier time $s = t - h$ form a Cox process $\zeta_s$ directed by $h^{-1}\xi_s$ (meaning that $\zeta_s$ is conditionally Poisson with intensity $h^{-1}\xi_s$, given $\xi_s$; cf. [13], page 226), and the generated clusters $\eta_h^i$ are conditionally independent and identically distributed apart from shifts. This is all explained in [8], pages 60ff, and some more precise statements with detailed proofs appear in Theorem 3.11 of [5] and Corollary 11.5.3 of [3]. In this paper, a generic cluster of age $t > 0$ is denoted by $\eta_t$; we write $P_x\{\eta_t \in \cdot\}$ for the distribution of a $t$-cluster centered at $x \in \mathbb{R}^d$ and put $P_\mu\{\eta_t \in \cdot\} = \int \mu(dx) P_x\{\eta_t \in \cdot\}$.

For the ease of reference, we state some basic scaling properties of DW-processes and their associated clusters (cf. Theorem 6.6 in [5]).

LEMMA 2.1. *Let $\xi$ be a DW-process in $\mathbb{R}^d$ starting at $\mu$ and with associated clusters $\eta_t$. Then for any $r, t > 0$, we have:*

(i) $\tilde{\xi}_t = r^{-2}\xi_{r^2 t} S_r$ *is a DW-process starting from $\tilde{\mu} = r^{-2}\mu S_r$,*
(ii) $\eta_t \stackrel{d}{=} r^{-2}\eta_{r^2 t} S_r$ *under $P_0$.*

PROOF. Part (i) may be proved by the argument in [8], page 51. A similar scaling property is then obtained for the cluster representation of $\xi$, and (ii) follows by the uniqueness of the associated Lévy measure (cf. Theorem 6.1 in [11]). □

Given a random measure $\xi$ on $\mathbb{R}^d$ with $\sigma$-finite intensity $E\xi$, we define the kernel of associated *Palm distributions* $Q_x$ by the disintegration formula

$$E \int f(x, \xi) \xi(dx) = \int E\xi(dx) \int f(x, \mu) Q_x(d\mu),$$

for any measurable function $f \geq 0$ on $\mathbb{R}^d \times \mathcal{M}_d$. If $\xi$ is defined on the canonical probability space with distribution $P$, we also write $P^x = Q_x$. When $\xi$ is stationary, we may choose the measures $P^0 = P^x \circ \theta_{-x}^{-1}$ to be independent of $x$, in which case $P^x = P^0 \circ \theta_x^{-1}$ for all $x$. What is said above applies even to the Palm distributions of pseudo-random measures $\tilde{\xi}$ on $\mathbb{R}^d$, as long as $\tilde{E}\tilde{\xi}$ is $\sigma$-finite. (In particular, the $\tilde{P}^x$ are still probability measures in this case, even if $\tilde{P}$ is not.)

In the nonstationary case, the Palm distributions $P^x$ are only determined for $x \in \mathbb{R}^d$ a.e. $E\xi$. However, the function $x \mapsto P^x$ may have a version with nice continuity properties. In Lemma 3.5 below, we show that when $\xi$ is a locally finite DW-process with initial measure $\mu$, the family of shifted Palm distributions $P_\mu^x \circ \theta_{-x}^{-1}$ can be chosen to be locally continuous in total variation. The continuous version is then unique, and the Palm distribution $P_\mu^0$ becomes well defined. This is the version with a nice probabilistic representation, given by Corollary 4.1.6 in [5] or Theorem 11.7.1 in [3].

In this paper, Palm distributions figure prominently only in Sections 8 and 9. The following uniform convergence criterion for shifted Palm measures will be needed in Section 8.

LEMMA 2.2. *Let $\xi$ and $\xi_n$ be random measures on $\mathbb{R}^d$ with locally finite intensities, where $\xi$ is stationary, and let $Q$ and $Q_s^n$ be versions of the associated shifted Palm distributions. Fix a set $B \in \hat{\mathcal{B}}^d$ satisfying:*

(i) $E\xi_n B \to E\xi B > 0$,
(ii) $\|E[\xi_n B; \xi_n \in \cdot] - E[\xi B; \xi \in \cdot]\| \to 0$,
(iii) $\sup_{r,s \in B} \|Q_r^n - Q_s^n\| \to 0$.

*Then $\sup_{s \in B} \|Q_s^n - Q\| \to 0$.*

PROOF. For measurable $A \subset \mathcal{M}_d$, we define
$$f_A(\mu) = (\mu B)^{-1} \int_B \mu(ds) 1_A(\mu\theta_s), \qquad \mu \in \mathcal{M}_d,$$
where $0^{-1} 0 = 0$. Then
$$\int_B E\xi_n(ds) Q_s^n A = E \int_B \xi_n(ds) 1_A(\xi_n\theta_s)$$
$$= E\xi_n B f_A(\xi_n) = \int E[\xi_n B; \xi_n \in d\mu] f_A(\mu),$$
and similarly for $\xi$ and $Q$. Since $|\nu f| \leq \|\nu\|$ for any signed measure $\nu$ and measurable function $f$ into $[0,1]$, we get for $s \in B$
$$E\xi B \|Q_s^n - Q\| \leq \|E\xi B Q_s^n - E\xi_n B Q_s^n\| + \left\|E\xi_n B Q_s^n - \int_B E\xi_n(dr) Q_r^n\right\|$$
$$+ \left\|\int_B E\xi_n(dr) Q_r^n - E\xi B Q\right\|$$
$$\leq |E\xi B - E\xi_n B| + \int_B E\xi_n(dr) \|Q_s^n - Q_r^n\|$$
$$+ \|E[\xi_n B; \xi_n \in \cdot] - E[\xi B; \xi \in \cdot]\|.$$
By (i)–(iii) the right-hand side tends to 0 as $n \to \infty$, uniformly in $s \in B$, and the assertion follows since $E\xi B > 0$. □

We conclude this section with an elementary but somewhat technical interpolation principle that will be needed in Section 7.

LEMMA 2.3. *Let the functions $f, g > 0$ on $(0,1]$ and constants $p, c > 0$ be such that $f$ is nondecreasing, $\log g(e^{-t})$ is bounded and uniformly continuous on $\mathbb{R}_+$, and $t^{-p} f(t) g(t) \to c$ as $t \to 0$ along every sequence $(r^n)$ with $r$ in some dense set $D \subset (0,1)$. Then the same convergence holds along $(0,1)$.*

PROOF. Letting $w$ be the modulus of continuity of $\log g(e^{-t})$, we get
$$e^{-w(h)}g(e^{-t}) \leq g(e^{-t-h}) \leq e^{w(h)}g(e^{-t}), \qquad t, h \geq 0.$$
Writing $b_r = \exp w(-\log r)$, we obtain
$$b_r^{-1} g(t) \leq g(rt) \leq b_r g(t), \qquad t, r \in (0, 1).$$
For any $r, t \in (0,1)$, define $n = n(r,t)$ by $r^{n+1} < t \leq r^n$. Then by the monotonicity of $f$
$$r^p(r^{n+1})^{-p} f(r^{n+1}) b_r^{-1} g(r^{n+1}) \leq t^{-p} f(t) g(t)$$
$$\leq r^{-p}(r^n)^{-p} f(r^n) b_r g(r^n).$$
Letting $t \to 0$ for fixed $r \in D$, we get by the hypothesis
$$r^p b_r^{-1} c \leq \liminf_{t \to 0} t^{-p} f(t) g(t) \leq \limsup_{t \to 0} t^{-p} f(t) g(t) \leq r^{-p} b_r c.$$
It remains to note that $r^{-p} b_r \to 1$ as $r \to 1$ along $D$. □

**3. Moments and continuity.** Throughout the paper we need some basic results involving the first and second moment measures $E_\mu \xi_t$ and $E_\mu \xi_t^2$ of a DW-process $\xi$ in $\mathbb{R}^d$. Here a simple estimate for the normal densities $p_t$ will be useful.

LEMMA 3.1. *Let $p_t$, $t > 0$, denote the standard normal density functions on $\mathbb{R}^d$. Then for fixed $d$ and $t$ we have*
$$p_t(x+y) \lesssim p_{t+h}(x), \qquad x \in \mathbb{R}^d, |y| \leq h \leq t.$$

PROOF. If $|x| \geq 4t$ and $|y| \leq h$, then $|y|/|x| \leq h/4t$, and so for $r = h/t \leq 1$
$$\frac{|x+y|^2}{t} \frac{(t+h)}{|x|^2} \geq \left(1 - \frac{|y|}{|x|}\right)^2 \left(1 + \frac{h}{t}\right) \geq \left(1 - \frac{r}{2}\right)(1+r) \geq 1,$$
which implies $p_t(x+y) \lesssim p_{t+h}(x)$ when $h \leq t$. The same relation holds trivially for $|x| \leq 4t$ and $|y| \leq h \leq t$. □

Let us now consider the intensity measures $E_\mu \xi_t$ of a DW-process $\xi$ starting from an arbitrary $\sigma$-finite measure $\mu$.

LEMMA 3.2. *Let $\xi$ be a DW-process in $\mathbb{R}^d$ with associated clusters $\eta_t$, and fix a $\sigma$-finite measure $\mu$. Then for any fixed $t > 0$, the following two conditions are equivalent:*

(i) *$\xi_t$ is locally finite a.s. $P_\mu$,*
(ii) *$E_\mu \xi_t$ is locally finite.*

*Furthermore*, (i) *and* (ii) *hold for all* $t > 0$ *iff*:

(iii) $\mu p_t < \infty$ *for all* $t > 0$,

*in which case we have for any* $t > 0$:

(iv) $E_\mu \xi_t = t^{-1} E_\mu \eta_t$ *has the finite, continuous density* $\mu * p_t$,
(v) $E_\mu(\xi_t \theta_x)$ *is locally continuous in total variation in* $x$, *and the same continuity holds globally when* $\mu$ *is bounded*.

PROOF. The formula $E_\mu \xi_t = (\mu * p_t) \cdot \lambda^d$, well known for bounded $\mu$ (cf. Lemma 2.1 in [8]), extends by monotone convergence to any $\sigma$-finite measure $\mu$ (though $E_\mu \xi_t$ may fail to be $\sigma$-finite, in general). The relation $E_\mu \eta_t = t E_\mu \xi_t$ follows from the cluster representation of $\xi_t$.

Condition (ii) clearly implies (i). Conversely, let $B = B_x^\varepsilon$ with $\varepsilon^2 < t$ and $0 < E_\mu \xi_t B < \infty$. Using the Paley–Zygmund inequality (cf. [13], page 63) and Hint (2) in [26], page 239, we get for any $r \in (0,1)$

$$P_\mu\left\{\frac{\xi_t B}{E_\mu \xi_t B} > r\right\} \geq (1-r)^2 \frac{(E_\mu \xi_t B)^2}{E_\mu(\xi_t B)^2} \geq \frac{(1-r)^2}{1 + c(E_\mu \xi_t B)^{-1}},$$

for some constant $c > 0$ depending only on $d$, $t$ and $\varepsilon$. Now assume instead that $E_\mu \xi_t B = \infty$, and choose some bounded measures $\mu_n \uparrow \mu$ with $E_{\mu_n} \xi_t B > n$. Applying the previous inequality to each $\mu_n$ gives

$$P_\mu\{\xi_t B > rn\} \geq P_{\mu_n}\{\xi_t B > r E_{\mu_n} \xi_t B\} \geq \frac{(1-r)^2}{1 + c/n}.$$

Letting $n \to \infty$ and then $r \to 0$, we obtain $\xi_t B = \infty$ a.s. In particular, this shows that (i) implies (ii).

Next assume (iii). Fixing $x \in \mathbb{R}^d$ and choosing $r \geq t + 2|x|$, we see from Lemma 3.1 that $p_t(x - u) \lesssim p_r(u)$ and hence $(\mu * p_t)(x) \lesssim \mu p_r < \infty$, which shows that $E_\mu \xi_t$ has the finite density $\mu * p_t$. Next we may write

$$|(\mu * p_t)(x + y) - (\mu * p_t)(x)| \leq \int |p_t(x + y - u) - p_t(x - u)| \mu(du),$$

where the integrand tends to 0 as $y \to 0$. Furthermore, Lemma 3.1 yields

(1) $\qquad |p_t(x + y - u) - p_t(x - u)| \lesssim p_{2t}(x - u), \qquad |y| \leq t.$

Since $\mu * p_{2t}(x) < \infty$, the continuity of $\mu * p_t$ follows by dominated convergence. This proves (iv), which in turn implies (ii) for every $t > 0$. Conversely, (ii) yields $(\mu * p_n)(x) < \infty$ for all $n \in \mathbb{N}$ and for $x \in \mathbb{R}^d$ a.e. $\lambda^d$. Fixing such an $x$ and using Lemma 3.1 as before, we obtain condition (iii).

To prove (v), we write for any $y \in \mathbb{R}^d$ and $t > 0$

$$\|E_\mu(\xi_t \theta_y) - E_\mu \xi_t\| = \int |(\mu * p_t)(x - y) - (\mu * p_t)(x)| \, dx$$

$$\leq \int \mu(du) \int |p_t(x - y - u) - p_t(x - u)| \, dx,$$

where the integrand tends to 0 as $y \to 0$. For bounded $\mu$, we may use (1) again and note that $\int \mu(du) \int p_{2t}(x-u)\,dx = \|\mu\| < \infty$, which justifies taking limits under the integral sign. For general $\mu$ as in (iii), fix any $B \in \hat{\mathcal{B}}^d$, and note that

$$\|E_\mu(\xi_t \theta_y) - E_\mu \xi_t\|_B \leq \int \mu(du) \int_B |p_t(x-y-u) - p_t(x-u)|\,dx.$$

Choosing $r > 0$ so large that $t + 2|x-y| \leq r$ for $x \in B$ and $|y| \leq 1$, we see from Lemma 3.1 that $p_t(x-y-u) \lesssim p_r(u)$ for any such $x$ and $y$. Since $\mu p_r < \infty$ by (iii), this justifies the dominated convergence in this case, and (v) follows. □

Assuming the DW-process $\xi$ in $\mathbb{R}^d$ to be *locally finite under* $P_\mu$, in the sense that condition (i) above holds for all $t > 0$, we go on to study the second moment measures $E_\mu \xi_t^2$ and the associated *covariance measures* $\text{Cov}_\mu \xi_t = E_\mu \xi_t^2 - (E_\mu \xi_t)^2$ on $\mathbb{R}^{2d}$.

LEMMA 3.3. *Let the DW-process $\xi$ in $\mathbb{R}^d$ be locally finite under $P_\mu$, and denote the associated clusters by $\eta_t$. Then for any $t > 0$ and $x_1, x_2 = \bar{x} \pm r$ in $\mathbb{R}^d$, we have:*

(i) $\text{Cov}_\mu \xi_t = t^{-1} E_\mu \eta_t^2 = (\mu * q_t) \cdot \lambda^{2d}$ with $q_t = 2\int_0^t (p_s * p_{t-s}^{\otimes 2})\,ds$,
(ii) $(\mu * q_t)(x_1, x_2) \lesssim (\mu * p_t)(\bar{x}) p_t(r) |r|^{2-d} t^{d/2}$ for $d \geq 3$,
(iii) $(\mu * q_t)(x_1, x_2) \lesssim (\mu * p_t)(\bar{x})(t p_{t/2}(r) + \log_+(t/|r|^2))$ for $d = 2$,
(iv) $(\mu * q_t)(x_1, x_2) \asymp (\mu * p_t)(x)|\log|r||$ for $d = 2$ as $x_1, x_2 \to x$.

Note that the convolutions in (i) are defined, for any $x, y \in \mathbb{R}^d$, by

$$(\mu * q_t)(x, y) = \int \mu(du) q_t(x-u, y-u),$$

$$(p_s * p_{t-s}^{\otimes 2})(x, y) = \int p_s(u) p_{t-s}(x-u) p_{t-s}(y-u)\,du.$$

PROOF. (i) The expression for $\text{Cov}_\mu \xi_t$, well known for bounded $\mu$ (cf. [8], page 37f), extends to the general case by monotone convergence. To see that $E_\mu \eta_t^2 = t \text{Cov}_\mu \xi_t$, let $\zeta_0$ be the process of ancestors of $\xi_t$ at time 0, and denote the generated clusters by $\eta_t^i$. Using the Poisson property of $\zeta_0$ and the conditional independence of the clusters, we get

$$E_\mu \xi_t^2 = (E_\mu \xi_t)^2 + \text{Cov}_\mu \xi_t = E_\mu \sum_{i,j}(\eta_t^i \otimes \eta_t^j)$$

$$= \iint_{x \neq y} E_\mu \zeta_0^2(dx\,dy)(E_x \eta_t \otimes E_y \eta_t) + \int E_\mu \zeta_0(dx) E_x \eta_t^2$$

$$= t^{-2}(E_\mu \eta_t)^2 + t^{-1} E_\mu \eta_t^2 = (E_\mu \xi_t)^2 + t^{-1} E_\mu \eta_t^2.$$

(ii) By definition

$$(2) \quad q_t(x_1, x_2) = 2 \int_0^t ds \int p_s(u) p_{t-s}(x_1 - u) p_{t-s}(x_2 - u) \, du.$$

To estimate $q_t$, we may use the parallelogram identity to get

$$p_t(x_1) p_t(x_2) \equiv t^{-d} \exp(-(|x_1|^2 + |x_2|^2)/2t)$$
$$= t^{-d} \exp(-(|\bar{x}|^2 + |r|^2)/t) \equiv p_{t/2}(\bar{x}) p_{t/2}(r).$$

Applying this to (2) and using the semigroup property of the normal densities, we obtain

$$q_t(x_1, x_2) \equiv \int_0^t ds \int p_s(u) p_{(t-s)/2}(\bar{x} - u) p_{(t-s)/2}(r) \, du$$
$$= \int_0^t p_{(t-s)/2}(r) p_{(t+s)/2}(\bar{x}) \, ds$$
$$\lesssim p_t(\bar{x}) \int_0^t p_{s/2}(r) \, ds \equiv p_t(\bar{x}) \int_0^t s^{-d/2} e^{-|r|^2/s} \, ds$$
$$= p_t(\bar{x}) |r|^{2-d} \int_{|r|^2/t}^\infty v^{d/2-2} e^{-v} \, dv$$
$$\lesssim p_t(\bar{x}) |r|^{2-d} e^{-|r|^2/2t} \equiv p_t(\bar{x}) p_t(r) |r|^{2-d} t^{d/2}.$$

The required estimate now follows by convolution with $\mu$.

(iii) Here we see as before that

$$q_t(x_1, x_2) \lesssim p_t(\bar{x}) \int_{|r|^2/t}^\infty v^{-1} e^{-v} \, dv.$$

For $|r|^2 \leq t/2$ we have

$$\int_{|r|^2/t}^\infty v^{-1} e^{-v} \, dv \lesssim \int_{|r|^2/t}^1 v^{-1} \, dv = \log(t/|r|^2),$$

and for $|r|^2 \geq t/2$ we get

$$\int_{|r|^2/t}^\infty v^{-1} e^{-v} \, dv \leq \int_{|r|^2/t}^\infty e^{-v} \, dv = \exp(-|r|^2/t) \lesssim t p_{t/2}(r).$$

(iv) For fixed $\varepsilon > 0$ we have

$$q_t(x_1, x_2) \equiv \int_0^t p_{(t-s)/2}(r) p_{(t+s)/2}(\bar{x}) \, ds$$
$$\sim \int_{t-\varepsilon}^t p_{(t-s)/2}(r) p_{(t+s)/2}(\bar{x}) \, ds,$$

since
$$\int_0^{t-\varepsilon} p_{(t-s)/2}(r) p_{(t+s)/2}(\bar{x}) \, ds \lesssim p_t(\bar{x}) \int_{|r|^2/t}^{|r|^2/\varepsilon} v^{-1} \, dv$$
$$\to p_t(x) \log(t/\varepsilon) < \infty.$$

Noting that $p_{(t+s)/2}(\bar{x}) \to p_t(\bar{x}) \to p_t(x)$ as $s \to t$ and then $x_1, x_2 \to x$, we get for fixed $b > 0$

$$q_t(x_1, x_2) \asymp p_t(x) \int_0^t p_{s/2}(r) \, ds = p_t(x) \int_{|r|^2/t}^\infty v^{-1} e^{-v} \, dv$$
$$\sim p_t(x) \int_{|r|^2/t}^b v^{-1} e^{-v} \, dv,$$

where the last relation holds since $\int_b^\infty v^{-1} e^{-v} \, dv < \infty$. Since $e^{-v} \to 1$ as $v \to 0$, we obtain

$$q_t(x_1, x_2) \asymp p_t(x) \int_{|r|^2/t}^1 v^{-1} \, dv = p_t(x) \log(t/|r|^2) \asymp p_t(x) |\log |r||.$$

This proves the assertion for $\mu = \delta_0$. For general $\mu$, let $c > 0$ be such that $q_t(x_1, x_2) \sim c p_t(x) |\log |r||$. We need to show that

$$|\log |r||^{-1} \int \mu(du) q_t(x_1 - u, x_2 - u) \to c \int \mu(du) p_t(x - u),$$

as $x_1, x_2 \to x$ for fixed $\mu$ and $t$. Then note that by (iii) and Lemma 3.1

$$|\log |r||^{-1} q_t(x_1 - u, x_2 - u) \lesssim p_t(\bar{x} - u) \lesssim p_{t+h}(x - u),$$

as long as $|\bar{x} - x| \leq h$. Since $(\mu * p_{t+h})(x) < \infty$, the desired relation follows by dominated convergence. □

Part (iv) of the last lemma yields a useful scaling property for the second moments of a DW-process in $\mathbb{R}^2$. This will be needed in Section 9.

LEMMA 3.4. *Let the DW-process $\xi$ in $\mathbb{R}^2$ be locally finite under $P_\mu$. Consider a measurable function $f \geq 0$ on $\mathbb{R}^4$ such that $f(x,y) \log(|x-y|^{-1} \vee e)$ is integrable, where $x, y \in \mathbb{R}^2$, and suppose that either $\mu$ or $\mathrm{supp}\, f$ is bounded. Then as $\varepsilon \to 0$ for fixed $t > 0$, we have*

$$E_\mu (\xi_t S_\varepsilon)^2 f \asymp \varepsilon^4 |\log \varepsilon| \lambda^4 f \mu p_t.$$

*This holds in particular when both $f$ and its support are bounded. The statement remains true with $\xi_t$ replaced by the associated clusters $\eta_t$.*

PROOF. By Lemma 3.3(iv), the density $g$ of $E_\mu \xi_t^2$ satisfies

$$g(x_1, x_2) \sim c|\log |x_1 - x_2||(\mu * p_t)(\tfrac{1}{2}(x_1 + x_2)), \qquad x_1 \approx x_2 \text{ in } \mathbb{R}^2,$$

for some constant $c > 0$, and is otherwise bounded for bounded $\mu$. Furthermore, we have

$$E_\mu \xi_t^2 f = \int f(u/\varepsilon) g(u)\, du = \varepsilon^4 |\log \varepsilon| \int f(x) \frac{g(\varepsilon x)}{|\log \varepsilon|}\, dx.$$

Here the ratio in the last integrand tends to $c\mu p_t$ as $\varepsilon \to 0$. If $\mu$ or $\operatorname{supp} f$ is bounded, then the integral tends to $c\mu p_t \lambda^4 f$ by dominated convergence. To check the stated integrability condition when $f$ is bounded, we may change $(x_1, x_2)$ into the new coordinates $x_1 \pm x_2$, then replace $x_1 - x_2$ by polar coordinates $(r, \theta)$ and note that $\int_0^1 r|\log r|\, dr < \infty$. □

Next we prove the strong continuity under shifts for the distributions of a DW-process and the associated Palm distributions. This result will be needed in Sections 8 and 9.

LEMMA 3.5. *Let the DW-process $\xi$ in $\mathbb{R}^d$ be locally finite under $P_\mu$. Then for fixed $t > 0$, the distributions $P_\mu\{\xi_t \theta_x \in \cdot\}$ and $P_\mu^x\{\xi_t \theta_x \in \cdot\}$ are continuous in total variation, locally in $x \in \mathbb{R}^d$. The continuity holds even globally when $\mu$ is bounded.*

PROOF. First let $\|\mu\| < \infty$. Let $\zeta_s$ denote the ancestral process at time $s \in [0, t]$, and put $\tau = \inf\{s > 0; \|\zeta_s\| > \|\zeta_0\|\}$. Then $\zeta_0$ is Poisson with intensity $\mu/t$, and each ancestor in $\zeta_0$ branches before time $h \in (0, t)$ with probability $h/t$. Hence, the number of such branching individuals is Poisson with mean $\|\mu\|h/t^2$, and so $P\{\tau > h\} = \exp(-\|\mu\|h/t^2)$. Conditionally on $\tau > h$, the process $\zeta_h$ is again Poisson with intensity $t^{-1}(\mu * p_h) \cdot \lambda^d = E_\mu \xi_h / t$, and $\xi_t$ is conditionally independent of the event $\{\tau > h\}$, given $\zeta_h$. Therefore,

$$\|P_\mu\{\xi_t \theta_r \in \cdot\} - P_\mu\{\xi_t \in \cdot\}\|$$
$$\leq P_\mu\{\tau \leq h\} + \|P_\mu[\zeta_h \theta_r \in \cdot | \tau > h] - P_\mu[\zeta_h \in \cdot | \tau > h]\|$$
$$\leq (1 - e^{-\|\mu\|h/t^2}) + t^{-1}\|E_\mu \xi_h \theta_r - E_\mu \xi_h\|,$$

which tends to 0 as $r \to 0$ and then $h \to 0$ by Lemma 3.2(v).

For general $\mu$, we may choose some bounded measures $\mu_n \uparrow \mu$, so that $\mu_n' = \mu - \mu_n \downarrow 0$. Fixing any $B \in \hat{\mathcal{B}}^d$, we have

$$\|P_\mu\{\xi_t \theta_r \in \cdot\} - P_\mu\{\xi_t \in \cdot\}\|_B$$
$$\leq \|P_{\mu_n}\{\xi_t \theta_r \in \cdot\} - P_{\mu_n}\{\xi_t \in \cdot\}\| + 2P_{\mu_n'}\{\xi_t(B \cup \theta_r B) > 0\},$$

which tends to 0 as $r \to 0$ and then $n \to \infty$, by the previous case and the simple Lemma 4.3 below (whose proof is independent of the present result). This yields the continuity of $P_\mu\{\xi_t \theta_r \in \cdot\}$.

We turn to the Palm distributions $P_\mu^0\{\xi_t \in \cdot\}$. By Lemma 10.6 in [11] (cf. Lemma 11.4.2 in [3]), the measure $P_\mu^0\{\xi_t \in \cdot\}$ is the convolution of $P_\mu\{\xi_t \in \cdot\}$ with the Palm distribution at 0 of the Lévy measure $P_\mu\{\eta_t \in \cdot\} = \int \mu(dx) P_x\{\eta_t \in \cdot\}$. By the previous result and Fubini's theorem, it is then enough to show that the latter factor is continuous in total variation under shifts in $\mu$. By Corollary 4.1.6 in [5] (cf. Theorem 11.7.1 in [3]), the corresponding historical path is a Brownian bridge $X$ on $[0,t]$ from $\alpha$ to 0, where $\alpha$ has distribution $(p_t \cdot \mu)/\mu p_t$. The measure $\eta_t$ is the sum of independent clusters rooted along the path of $X$, with birth times given by an independent Poisson process $\zeta$ on $[0,t]$ with rate $2/(t-s)$ at time $s$.

Let $\tau$ be the first point of $\zeta$. Since $P\{\tau \leq h\} \to 0$ as $h \to 0$ and since the event $\tau > h$ is independent of the restriction of $\zeta$ to the interval $[h,t]$, it suffices, for any fixed $h > 0$, to prove the continuity in total variation for the sum of clusters born after time $h$. Since $X$ is again a Brownian bridge on $[h,t]$, conditionally on $\alpha$ and $X_h$, the mentioned sum is conditionally independent of $\alpha$ given $X_h$, and it is enough to prove that $P_\mu\{X_h \in \cdot\}$ is continuous in total variation under shifts in $\mu$.

Then put $s = t - h$, and note that $X_h$ is conditionally $N(s\alpha, sh)$ given $\alpha = X_0$. Thus, the conditional density of $X_h$ equals $p_{sh}(x - s\alpha)$. Since $\alpha$ has density $(p_t \cdot \mu)/\mu p_t$, the unconditional density of $X_h$ becomes

$$f_\mu(x) = (\mu p_t)^{-1} \int p_{sh}(x - su) p_t(u) \mu(du), \qquad x \in \mathbb{R}^d.$$

Replacing $\mu$ by the shifted measure $\mu\theta_r$ yields the density

$$f_{\mu\theta_r}(x) = ((\mu * p_t)(r))^{-1} \int p_{sh}(x - su + sr) p_t(u - r) \mu(du),$$

and we need to show that $f_{\mu\theta_r} \to f_\mu$ in $L^1$ as $r \to 0$. Since $(\mu * p_t)(r) \to \mu p_t$ by Lemma 3.2(iv), it is enough to prove convergence of the $\mu$-integrals. Here the $L^1$-distance is bounded by

$$\int dx \int \mu(du) |p_{sh}(x - su + sr) p_t(u - r) - p_{sh}(x - su) p_t(u)|,$$

which tends to 0 as $r \to 0$ by Lemma 1.32 in [13], since the integrand tends to 0 by continuity and

$$\int dx \int \mu(du) p_{sh}(x - su + sr) p_t(u - r)$$
$$= (\mu * p_t)(r) \to \mu p_t = \int dx \int \mu(du) p_{sh}(x - su) p_t(u),$$

by Fubini's theorem and Lemma 3.2(iv). $\square$

**4. Hitting bounds and extinction.** In this section we derive some hitting estimates at fixed times for a DW-process $\xi$ in $\mathbb{R}^d$ and the associated clusters $\eta_t$. Those results will be useful throughout the remainder of the paper. We also discuss some extinction and related properties for DW-processes of dimension $d \geq 2$. For the ease of reference, we begin with a well-known relationship between the hitting probabilities of $\xi_t$ and $\eta_t$. Here and below $P_\mu\{\eta_t \in \cdot\} = \int \mu(dx) P_x\{\eta_t \in \cdot\}$.

LEMMA 4.1. *Let the DW-process $\xi$ in $\mathbb{R}^d$ with associated clusters $\eta_t$ be locally finite under $P_\mu$, and fix any $B \in \mathcal{B}^d$. Then*

$$P_\mu\{\eta_t B > 0\} = -t \log(1 - P_\mu\{\xi_t B > 0\}),$$
$$P_\mu\{\xi_t B > 0\} = 1 - \exp(-t^{-1} P_\mu\{\eta_t B > 0\}).$$

*In particular, $P_\mu\{\xi_t B > 0\} \sim t^{-1} P_\mu\{\eta_t B > 0\}$ as either side tends to 0.*

PROOF. Under $P_\mu$ we have $\xi_t = \sum_i \eta_t^i$, where the $\eta_t^i$ are conditionally independent clusters of age $t$ rooted at the points of a Poisson process with intensity $\mu/t$. For a cluster rooted at $x$, the hitting probability is $b_x = P_x\{\eta_t B > 0\}$. Hence (e.g., by Proposition 12.3 in [13]), the number of clusters hitting $B$ is Poisson distributed with mean $\mu b/t$, and so $P_\mu\{\xi_t B = 0\} = \exp(-\mu b/t)$, which yields the asserted formulas. □

Next we extend some classical hitting estimates for DW-processes of dimension $d \geq 2$. By Lemma 4.1 it is enough to consider the corresponding clusters $\eta_t$, and by shifting it suffices to consider balls centered at the origin.

LEMMA 4.2. *Let the $\eta_t$ be clusters of a DW-process in $\mathbb{R}^d$, and consider a $\sigma$-finite measure $\mu$ on $\mathbb{R}^d$.*

(i) *For $d \geq 3$, let $t_\varepsilon = t + \varepsilon^2$. Then for $0 < \varepsilon \leq \sqrt{t}$, we have*

$$\mu p_t \lesssim t^{-1} \varepsilon^{2-d} P_\mu\{\eta_t B_0^\varepsilon > 0\} \lesssim \mu p_{t(\varepsilon)}.$$

(ii) *For $d = 2$, we may choose $0 \leq l_\varepsilon - 1 \lesssim |\log \varepsilon|^{-1/2}$ and put $t_\varepsilon = t l_{\varepsilon/\sqrt{t}}$, so that uniformly for $x \in \mathbb{R}^2$ and $0 < \varepsilon < \frac{1}{2}\sqrt{t}$*

$$\mu p_t \lesssim t^{-1} \log(t/\varepsilon^2) P_\mu\{\eta_t B_0^\varepsilon > 0\} \lesssim \mu p_{t(\varepsilon)}.$$

PROOF. (i) For bounded $\mu$ we have by Theorem 3.1 in [4] (cf. Theorem III.5.11 and Exercise III.5.2 in [26])

$$\mu p_t \lesssim \varepsilon^{2-d} P_\mu\{\xi_t B_0^\varepsilon > 0\} \lesssim \mu p_{t(\varepsilon)},$$

and the asserted relations follow by Lemma 4.1. The result extends by linearity to any $\sigma$-finite measure $\mu$.

(ii) It is enough to take $t=1$, since by Lemma 2.1(ii) we then obtain for general $t>0$

$$\begin{aligned}P_x\{\eta_t B_0^\varepsilon > 0\} &= P_0\{\eta_t B_x^\varepsilon > 0\} = P_0\{\eta_1 B_{x/\sqrt{t}}^{\varepsilon/\sqrt{t}} > 0\}\\&\lesssim |\log(\varepsilon/\sqrt{t})|^{-1} p_{l(\varepsilon/\sqrt{t})}(x/\sqrt{t})\\&\lesssim t(\log(t/\varepsilon^2))^{-1} p_{t(\varepsilon)}(x),\end{aligned}$$

and similarly for the lower bound.

For $t=1$ we have by Theorem 2 in [21]

$$\begin{aligned}p_1(x) &\lesssim |\log\varepsilon| P_x\{\eta_1 B_0^\varepsilon > 0\}\\&\lesssim (1 + 1\{|x|^2 > |\log\varepsilon|\}|x|^4) p_1(|x|-\varepsilon).\end{aligned}$$

In particular, this gives the required lower bound. Next, Lemma 3.1 yields $p_1(|x|-\varepsilon) \lesssim p_{1+\varepsilon}(x)$, and by elementary estimates we get for $|\log\varepsilon| \geq e$

$$1 + 1\{|x|^2 > |\log\varepsilon|\}|x|^4 \lesssim \exp\left(\frac{2\log|\log\varepsilon|}{|\log\varepsilon|}|x|^2\right).$$

Hence, by combination, we get for $\varepsilon$ bounded by some constant $c>0$

$$|\log\varepsilon| P_x\{\eta_1 B_0^\varepsilon > 0\} \lesssim \exp\left\{-\frac{|x|^2}{2}\left(\frac{1}{1+\varepsilon} - \frac{4\log|\log\varepsilon|}{|\log\varepsilon|}\right)\right\} \lesssim p_{l(\varepsilon)}(x),$$

where

$$l(\varepsilon) = \left(\frac{1}{1+\varepsilon} - \frac{4\log|\log\varepsilon|}{|\log\varepsilon|}\right)^{-1}, \qquad 0 < \varepsilon \leq c.$$

As $\varepsilon \to 0$, we note that

$$0 \leq l(\varepsilon) - 1 \lesssim \varepsilon + \frac{4\log|\log\varepsilon|}{|\log\varepsilon|} \lesssim |\log\varepsilon|^{-1/2}.$$

When $c < \varepsilon < \frac{1}{2}$, we have instead

$$\begin{aligned}|\log\varepsilon| P_x\{\eta_1 B_0^\varepsilon > 0\} &\lesssim (1+|x|^4) p_1(|x|-\varepsilon)\\&\lesssim \exp(a|x|^2) p_{1+\varepsilon}(x),\end{aligned}$$

for any fixed $a > 0$. Choosing $a$ small enough, we get again a bound of the form $p_{l(\varepsilon)}$, for a suitable choice of $l(\varepsilon) \geq 1$. $\square$

The following simple result is often useful to extend results for bounded initial measures $\mu$ to the general case.

LEMMA 4.3. *Let the DW-process $\xi$ in $\mathbb{R}^d$ be locally finite under $P_\mu$, and suppose that $\mu \geq \mu_n \downarrow 0$. Then $P_{\mu_n}\{\xi_t B > 0\} \to 0$ as $n \to \infty$ for fixed $t > 0$ and $B \in \hat{\mathcal{B}}^d$.*

PROOF. We may assume that $B = B_0^r$ for some $r > 0$. Using Lemmas 3.2, 4.1 and 4.2, along with a projection argument when $d = 1$, we get for small enough $\varepsilon > 0$ and for suitable $t_\varepsilon > 0$

$$P_\mu\{\xi_t B > 0\} \lesssim \int_B P_\mu\{\xi_t B_x^\varepsilon > 0\}\, dx \lesssim \int_B (\mu * p_{t(\varepsilon)})(x)\, dx < \infty.$$

The assertion now follows by dominated convergence. □

Next we need to estimate the probability that a small ball in $\mathbb{R}^d$ is hit by more than one subcluster of our DW-process $\xi$. This result will play a crucial role throughout the remainder of the paper.

LEMMA 4.4. *Let the DW-process $\xi$ in $\mathbb{R}^d$ be locally finite under $P_\mu$. For any $t \geq h > 0$ and $\varepsilon > 0$, let $\kappa_h^\varepsilon$ be the number of $h$-clusters hitting $B_0^\varepsilon$ at time $t$. Then:*

(i) *for $d \geq 3$ and as $\varepsilon^2 \ll h \leq t$, we have with $t_\varepsilon = t + \varepsilon^2$*

$$E_\mu \kappa_h^\varepsilon (\kappa_h^\varepsilon - 1) \lesssim \varepsilon^{2(d-2)}(h^{1-d/2}\mu p_t + (\mu p_{t(\varepsilon)})^2),$$

(ii) *for $d = 2$ we may choose $0 < t_{h,\varepsilon} - t \lesssim h|\log \varepsilon|^{-1/2}$, such that as $\varepsilon \ll h \leq t$*

$$E_\mu \kappa_h^\varepsilon (\kappa_h^\varepsilon - 1) \lesssim |\log \varepsilon|^{-2}(\log(t/h)\mu p_t + (\mu p_{t(h,\varepsilon)})^2).$$

PROOF. (i) Let $\zeta_s$ be the Cox process of ancestors to $\xi_t$ at time $s = t - h$, and write $\eta_h^i$ for the associated $h$-clusters. Using Lemma 4.2(i), the conditional independence of the clusters and the fact that $E_\mu \zeta_s^2 = h^{-2} E_\mu \xi_s^2$ outside the diagonal, we get with $p_h^\varepsilon(x) = P_x\{\eta_h B_0^\varepsilon > 0\}$

$$E_\mu \kappa_h^\varepsilon(\kappa_h^\varepsilon - 1) = E_\mu \sum_{i \neq j} 1\{\eta_h^i B_0^\varepsilon \wedge \eta_h^j B_0^\varepsilon > 0\}$$

$$= \iint_{x \neq y} p_h^\varepsilon(x) p_h^\varepsilon(y) E_\mu \zeta_s^2(dx\, dy)$$

$$\lesssim \varepsilon^{2(d-2)} \iint p_{h(\varepsilon)}(x) p_{h(\varepsilon)}(y) E_\mu \xi_s^2(dx\, dy).$$

By Lemma 3.2, Fubini's theorem and the semigroup property of $(p_t)$, we get

$$\int p_{h(\varepsilon)}(x) E_\mu \xi_s(dx) = \int p_{h(\varepsilon)}(x)(\mu * p_s)(x)\, dx$$

$$= \int \mu(du)(p_{h(\varepsilon)} * p_s)(u) = \mu p_{t(\varepsilon)}.$$

Next, we get by Lemma 3.3(i), Fubini's theorem, the properties of $(p_t)$ and the relations $t \leq t_\varepsilon \leq 2t - s$

$$\iint p_{h(\varepsilon)}(x) p_{h(\varepsilon)}(y) \operatorname{Cov}_\mu \xi_s(dx\,dy)$$
$$= 2 \iint p_{h(\varepsilon)}(x) p_{h(\varepsilon)}(y)\,dx\,dy \int \mu(du) \int_0^s dr$$
$$\times \int p_r(v-u) p_{s-r}(x-v) p_{s-r}(y-v)\,dv$$
$$= 2 \int \mu(du) \int_0^s dr \int p_r(u-v)(p_{t(\varepsilon)-r}(v))^2\,dv$$
$$\lesssim \int \mu(du) \int_0^s (t-r)^{-d/2}(p_r * p_{(t(\varepsilon)-r)/2})(u)\,dr$$
$$= \int \mu(du) \int_0^s (t-r)^{-d/2} p_{(t(\varepsilon)+r)/2}(u)\,dr$$
$$\lesssim \int p_t(u)\mu(du) \int_h^t r^{-d/2}\,dr \lesssim \mu p_t h^{1-d/2}.$$

The assertion follows by combination of these estimates.

(ii) Here we may proceed as before, with the following changes: Using Lemma 4.2(ii) instead of (i), we see that the factor $\varepsilon^{2(d-2)}$ should be replaced by $|\log(\varepsilon/\sqrt{h})|^{-2} \lesssim |\log \varepsilon|^{-2}$. In the last computation, we have now $\int_h^t r^{-1}\,dr = \log(t/h)$. Since $h_\varepsilon = hl_{\varepsilon/\sqrt{h}}$ with $0 \leq l_\varepsilon - 1 \lesssim |\log \varepsilon|^{-1/2}$, we may choose $t_{h,\varepsilon} = t + (h_\varepsilon - h)$ in the second term on the right. As for the estimates leading up to the first term, we note that the bound $t_{h,\varepsilon} + s \leq 2t$ remains valid for sufficiently small $\varepsilon/h$. □

Using the bounds in Lemma 4.2, we may improve some known extinction criteria for DW-processes of dimension $d \geq 2$.

THEOREM 4.5. *Let $\xi$ be a locally finite DW-process in $\mathbb{R}^d$, $d \geq 2$, with arbitrary initial distribution. Then these conditions are equivalent as $t \to \infty$:*

(i) $\xi_t \xrightarrow{d} 0$,
(ii) $\operatorname{supp} \xi_t \xrightarrow{d} \varnothing$,
(iii) $\begin{cases} \xi_0 p_t \xrightarrow{P} 0, & d \geq 3, \\ (\log t)^{-1} \xi_0 p_t \xrightarrow{P} 0, & d = 2. \end{cases}$

Already Dawson [2] noted that $\xi_t \xrightarrow{d} 0$ for a DW-process in $\mathbb{R}^2$ with $\xi_0 = \lambda^2$. The equivalence of (i) and (iii) was proved for $d = 2$ by Bramson, Cox and Greven [1] (see also [19]). Condition (ii) means that $1\{\xi B > 0\} \xrightarrow{P} 0$ for all $B \in \hat{\mathcal{B}}^2$. The corresponding a.s. convergence fails for $d = 2$ and $\xi_0 = \lambda^2$, for example, by the ergodic theorem in [9] (cf. Theorem 2.25 of [8]). However, for $d = 1$ such an a.s. result was obtained by Iscoe [10].

PROOF. First let $d = 2$. Using Lemmas 2.1(i), 4.1 and 4.2(ii), along with the properties of $p_t$, we get for any measure $\mu$ and constants $r, t, \varepsilon > 0$ with $t\varepsilon^2 = 1$

$$\begin{aligned} P_\mu\{\xi_t B_0^r > 0\} &= P_{\varepsilon^2 \mu S_{1/\varepsilon}}\{\xi_{\varepsilon^2 t} B_0^{r\varepsilon} > 0\} \\ &\lesssim \varepsilon^2 \mu S_{1/\varepsilon} p_{l(r\varepsilon)} |\log(r\varepsilon)|^{-1} \\ &\lesssim \varepsilon^2 |\log \varepsilon|^{-1} \mu(p_2 \circ S_\varepsilon) \lesssim (\log 2t)^{-1} \mu p_{2t}, \end{aligned}$$

since $1 \leq l_\varepsilon \leq 2$ for sufficiently small $\varepsilon > 0$. Combining with the corresponding lower bound gives

$$(\log t)^{-1} \mu p_t \wedge 1 \lesssim P_\mu\{\xi_t B_0^r > 0\} \lesssim (\log 2t)^{-1} \mu p_{2t} \wedge 1,$$

and so for a general initial distribution

$$E[(\log t)^{-1} \xi_0 p_t \wedge 1] \lesssim P\{\xi_t B_0^r > 0\} \lesssim E[(\log 2t)^{-1} \xi_0 p_{2t} \wedge 1].$$

As $t \to \infty$, we obtain $1\{\xi_t B_0^r > 0\} \xrightarrow{P} 0$ iff $(\log t)^{-1} \xi_0 p_t \xrightarrow{P} 0$, and the equivalence of (ii) and (iii) follows since $r$ was arbitrary.

For $d \geq 3$, we may use Lemma 4.2(i) instead to write

$$\begin{aligned} P_\mu\{\xi_t B_0^r > 0\} &= P_{\varepsilon^2 \mu S_{1/\varepsilon}}\{\xi_{\varepsilon^2 t} B_0^{r\varepsilon} > 0\} \\ &\lesssim \varepsilon^2 \mu S_{1/\varepsilon} p_{l(r\varepsilon)}(r\varepsilon)^{d-2} \lesssim \varepsilon^d t^{d/2} \mu p_{2t} = \mu p_{2t}, \end{aligned}$$

and similarly for the lower bound. Hence, for a general initial distribution,

$$E[\xi_0 p_t \wedge 1] \lesssim P\{\xi_t B_0^\varepsilon > 0\} \lesssim E[\xi_0 p_{2t} \wedge 1],$$

which shows again that (ii) and (iii) are equivalent. Since clearly (ii) implies (i), it remains to prove that (i) implies (iii).

Then put $B = B_0^1$, and suppose that $\xi$ is locally finite under $P_\mu$. Noting that $\text{Var}_\mu \xi_t B \lesssim E_\mu \xi_t B$ for $d \geq 3$ by Hint (2) in [26], page 239, we see as in the proof of Lemma 3.2 that

$$P_\mu\left\{\frac{\xi_t B}{E_\mu \xi_t B} > r\right\} \geq \frac{(1-r)^2}{1 + c(E_\mu \xi_t B)^{-1}}, \qquad r \in (0, 1),$$

where the constant $c > 0$ depends only on $d$. Hence, if $\xi_t B \xrightarrow{P} 0$ along some sequence $t_n \to \infty$, we get $E_\mu \xi_t B \to 0$ along the same sequence. Noting that $E_\mu \xi_t B \gtrsim \mu p_{t-1}$ by Lemma 3.1, we obtain $\mu p_{t_n-1} \to 0$.

For general $\xi_0$, (i) implies $\xi_t B \xrightarrow{P} 0$. Hence, for any $t_n \to \infty$ we have $\xi_t B \to 0$ a.s. along some subsequence $(t_{n'})$. Since this remains conditionally true given $\xi_0$, we see as before that $\xi_0 p_t \to 0$ a.s. along the shifted sequence $(t_{n'} - 1)$. Since the sequence $(t_n - 1)$ was arbitrary, $\xi_0 p_t \xrightarrow{P} 0$ follows by Lemma 4.2 in [13]. $\square$

In the stationary case, we can also estimate the rate of clustering. For a stationary random measure $\zeta$ on $\mathbb{R}^d$, the associated *sample intensity* $\bar{\zeta}$ is defined by $\bar{\zeta} \cdot \lambda^d = E[\zeta | \mathcal{I}]$, where $\mathcal{I}$ denotes the invariant $\sigma$-field.

PROPOSITION 4.6. *Let $\xi$ be a DW-process in $\mathbb{R}^2$, starting from a stationary random measure $\xi_0 \neq 0$ with sample intensity $\bar{\xi}_0 < \infty$ a.s. Then $P\{\xi_t B_0^r > 0\} \to 0$ as $t \to \infty$ iff $r^2/t \to 0$.*

PROOF. Letting $t\varepsilon^2 = 1$ and $r^2/t \to 0$, we get as in the previous proof
$$P_\mu\{\xi_t B_0^r > 0\} \lesssim \varepsilon^2 |\log(r\varepsilon)|^{-1} \mu(p_2 \circ S_\varepsilon)$$
$$\lesssim (\log(t/r^2))^{-1} \mu p_{2t}.$$
Hence, for a general initial distribution
$$P\{\xi_t B_0^r > 0\} \lesssim E[(\log(t/r^2))^{-1} \xi_0 p_{2t} \wedge 1],$$
which tends to 0 as $r^2 \ll t \to \infty$, since $\xi_0 p_{2t} \to \bar{\xi}_0 < \infty$ a.s. by Corollary 10.19 in [13].

Conversely, truncating $r\varepsilon$ at $\frac{1}{2}$, we get as before
$$P\{\xi_t B_0^r > 0\} \gtrsim E[|\log(r\varepsilon \wedge \tfrac{1}{2})|^{-1} \xi_0 p_t \wedge 1],$$
and so $P\{\xi_t B_0^r > 0\} \to 0$ implies
$$|\log(r\varepsilon \wedge \tfrac{1}{2})|^{-1} \bar{\xi}_0 \lesssim |\log(r\varepsilon \wedge \tfrac{1}{2})|^{-1} \xi_0 p_t \xrightarrow{P} 0.$$
Since $P\{\bar{\xi}_0 > 0\} > 0$, we get $|\log(r\varepsilon \wedge \tfrac{1}{2})| \to \infty$ and therefore $r^2/t \to 0$. $\square$

**5. Hitting asymptotics.** For a DW-process $\xi$ of dimension $d \geq 3$, we know from Theorem 3.1 of Dawson, Iscoe and Perkins [4] (cf. Remark III.5.12 in [26]) that, as $\varepsilon \to 0$ for fixed $t > 0$, $x \in \mathbb{R}^d$ and bounded $\mu$,

(3) $$\varepsilon^{2-d} P_\mu\{\xi_t B_x^\varepsilon > 0\} \to c_d (\mu * p_t)(x),$$

where $c_d > 0$ is a constant depending only on $d$, and the convergence is uniform for $x \in \mathbb{R}^d$ and for bounded $t^{-1}$ and $\|\mu\|$. Here we prove a similar result for $d = 2$, with $c_d$ replaced by a suitable normalizing function $m$.

Writing $p_h^\varepsilon(x) = P_x\{\eta_h B_0^\varepsilon > 0\}$, where $\eta_h$ denotes an $h$-cluster associated with a DW-process in $\mathbb{R}^d$, we define our normalizing function for $d = 2$ by
$$m(\varepsilon) = |\log \varepsilon| \lambda^2 p_1^\varepsilon = |\log \varepsilon| P_{\lambda^2}\{\eta_1 B_0^\varepsilon > 0\}, \qquad \varepsilon > 0.$$
The following technical result will play a crucial role below, especially in Section 7.

LEMMA 5.1. *The function $t \mapsto \log m(\exp(-e^t))$ is bounded and uniformly continuous on $[1, \infty)$.*

PROOF. The boundedness of $\log m$ is clear from Lemma 4.2(ii). For any $h \in (0,1]$, let $\zeta_s$ be the process of ancestors to $\xi_1$ at time $s = 1 - h$, and denote the generated $h$-clusters by $\eta_h^i$. Then for $0 < r \ll 1$ and $0 < \varepsilon \ll h$ we get the following chain of relations, to be explained and justified below:

$$m(\varepsilon)|\log \varepsilon|^{-1} \approx r^{-1} P_{r\lambda^2}\{\xi_1 B_0^\varepsilon > 0\}$$
$$\approx r^{-1} E_{r\lambda^2} \sum_i 1\{\eta_h^i B_0^\varepsilon > 0\} = r^{-1} E_{r\lambda^2} \zeta_s p_h^\varepsilon$$
$$= h^{-1} P_{\lambda^2}\{\eta_h B_0^\varepsilon > 0\} = P_{\lambda^2}\{\eta_1 B_0^{\varepsilon/\sqrt{h}} > 0\}$$
$$= m(\varepsilon/\sqrt{h})|\log(\varepsilon/\sqrt{h})|^{-1} \approx m(\varepsilon/\sqrt{h})|\log \varepsilon|^{-1}.$$

Here the first two steps are suggested by Lemmas 4.1 and 4.4(ii), respectively, the third step holds by the conditional independence of the clusters, the fourth step holds by the Cox property of $\zeta_s$, the fifth step holds by Lemma 2.1(ii), the sixth step holds by the definition of $m$ and the last step is suggested by the relation $\varepsilon \ll h$.

To estimate the approximation errors, we see from Lemmas 4.1 and 4.2(ii) that

$$|m(\varepsilon) - r^{-1}|\log \varepsilon| P_{r\lambda^2}\{\xi_1 B_0^\varepsilon > 0\}|$$
$$= r^{-1}|\log \varepsilon||P_{r\lambda^2}\{\eta_1 B_0^\varepsilon > 0\} - P_{r\lambda^2}\{\xi_1 B_0^\varepsilon > 0\}|$$
$$\lesssim r^{-1}|\log \varepsilon|(P_{r\lambda^2}\{\eta_1 B_0^\varepsilon > 0\})^2$$
$$\lesssim r^{-1}|\log \varepsilon|^{-1}(r\lambda^2 p_{l(\varepsilon)})^2 = r|\log \varepsilon|^{-1}.$$

Next, Lemma 4.4(ii) yields

$$r^{-1}|\log \varepsilon|\left| E_{r\lambda^2} \sum_i 1\{\eta_h^i B_0^\varepsilon > 0\} - P_{r\lambda^2}\{\xi_1 B_0^\varepsilon > 0\}\right|$$
$$= r^{-1}|\log \varepsilon| E_{r\lambda^2}(\kappa_h^\varepsilon - 1)_+$$
$$\lesssim \frac{|\log h|r\lambda^2 p_1 + (r\lambda^2 p_{t(h,\varepsilon)})^2}{r|\log \varepsilon|} = \frac{|\log h| + r}{|\log \varepsilon|}.$$

Finally, we note that

$$m(\varepsilon/\sqrt{h})\left|\frac{|\log \varepsilon|}{|\log(\varepsilon/\sqrt{h})|} - 1\right| \lesssim \frac{|\log h|}{|\log \varepsilon|},$$

by the boundedness of $m$. Combining these estimates and letting $r \to 0$, we obtain
$$|m(\varepsilon) - m(\varepsilon/\sqrt{h})| \lesssim \frac{|\log h|}{|\log \varepsilon|}.$$

Taking $\varepsilon = e^{-t}$ and $\varepsilon/\sqrt{h} = e^{-s}$ with $t - s \ll t$ gives
$$\left|\log \frac{m(e^{-t})}{m(e^{-s})}\right| \lesssim \left|\frac{m(e^{-t})}{m(e^{-s})} - 1\right| \lesssim |m(e^{-t}) - m(e^{-s})|$$
$$\lesssim (t-s)/t \lesssim |\log(t/s)|,$$

which extends immediately to arbitrary $s, t \geq 1$. Replacing $s$ and $t$ by $e^s$ and $e^t$ gives
$$|\log m(\exp(-e^t)) - \log m(\exp(-e^s))| \lesssim |t - s|,$$

which implies the asserted uniform continuity. $\square$

We proceed to approximate the hitting probabilities $p_h^\varepsilon$ by suitably normalized Dirac functions. Even this result will play a crucial role in the sequel, both here and in Section 7.

LEMMA 5.2. *Write $p_h^\varepsilon(x) = P_x\{\eta_h B_0^\varepsilon > 0\}$, where the $\eta_h$ are clusters of a DW-process in $\mathbb{R}^d$, and fix a bounded, uniformly continuous function $f \geq 0$ on $\mathbb{R}^d$. Then:*

(i) *for $d \geq 3$ and as $0 < \varepsilon^2 \ll h \to 0$, we have*
$$\|h^{-1}\varepsilon^{2-d}(p_h^\varepsilon * f) - c_d f\| \to 0,$$

(ii) *for $d = 2$ and as $0 < \varepsilon \leq h \to 0$ with $|\log h| \ll |\log \varepsilon|$, we have*
$$\|h^{-1}|\log \varepsilon|(p_h^\varepsilon * f) - m(\varepsilon) f\| \to 0.$$

*Both results hold uniformly over any class of uniformly bounded and equicontinuous functions $f \geq 0$ on $\mathbb{R}^d$.*

PROOF. (i) Using (3) and Lemmas 2.1(ii), 4.1 and 4.2(i), we get by dominated convergence
$$(4) \qquad \lambda^d p_h^\varepsilon = h^{d/2} \lambda^d p_1^{\varepsilon/\sqrt{h}} \sim c_d h^{d/2} (\varepsilon/\sqrt{h})^{d-2} \lambda^d p_1 = c_d \varepsilon^{d-2} h.$$

Similarly, Lemma 4.2(ii) yields for fixed $r > 0$ and a standard normal random vector $\gamma$ in $\mathbb{R}^d$
$$\varepsilon^{2-d} h^{-1} \int_{|x|>r} p_h^\varepsilon(x)\, dx \lesssim \int_{|u|>r/\sqrt{h}} p_{l(\varepsilon)}(u)\, du$$
$$(5) \qquad\qquad\qquad\qquad\qquad = P\{|\gamma| l_\varepsilon^{1/2} > r/\sqrt{h}\} \to 0.$$

By (4) it is enough to show that $\|\hat{p}_h^\varepsilon * f - f\| \to 0$ as $h$, $\varepsilon^2/h \to 0$, where $\hat{p}_h^\varepsilon = p_h^\varepsilon/\lambda^d p_h^\varepsilon$. Writing $w_f$ for the modulus of continuity of $f$, we get

$$\|\hat{p}_h^\varepsilon * f - f\| = \sup_x \left| \int \hat{p}_h^\varepsilon(u)(f(x-u) - f(x))\, du \right|$$

$$\leq \int \hat{p}_h^\varepsilon(u) w_f(|u|)\, du$$

$$\leq w_f(r) + 2\|f\| \int_{|u|>r} \hat{p}_h^\varepsilon(u)\, du,$$

which tends to 0 as $h$, $\varepsilon^2/h \to 0$ and then $r \to 0$, by (5) and the uniform continuity of $f$.

(ii) By Lemmas 2.1(ii) and 5.1 we have

$$\lambda^2 p_h^\varepsilon = h\lambda^2 p_1^{\varepsilon/\sqrt{h}} = hm(\varepsilon/\sqrt{h})|\log(\varepsilon/\sqrt{h})|^{-1} \sim hm(\varepsilon)|\log \varepsilon|^{-1}.$$

We also see that, with $t_\varepsilon$ as in Lemma 4.2(ii),

$$h^{-1}|\log \varepsilon| \int_{|x|>r} p_h^\varepsilon(x)\, dx \lesssim \int_{|u|>r/\sqrt{h}} p_{t(\varepsilon)}(u)\, du \to 0.$$

The proof may now be completed as in case of (i). The last assertion is clear from the estimates in the preceding proofs. $\square$

We may now prove the mentioned convergence of suitably normalized hitting probabilities, a result that is often needed in subsequent sections. The case $d \geq 3$ is included for convenience of reference.

THEOREM 5.3. *Let $\xi$ be a DW-process in $\mathbb{R}^d$. Then for any $t > 0$ and bounded $\mu$, we have as $\varepsilon \to 0$:*

(i) $\|\varepsilon^{2-d} P_\mu\{\xi_t B_\cdot^\varepsilon > 0\} - c_d(\mu * p_t)\| \to 0$ *for $d \geq 3$,*
(ii) $\||\log \varepsilon| P_\mu\{\xi_t B_\cdot^\varepsilon > 0\} - m(\varepsilon)(\mu * p_t)\| \to 0$ *for $d = 2$,*

*and similarly for the clusters $\eta_t$ with $p_t$ replaced by $tp_t$. The results hold locally whenever $\xi$ is locally finite under $P_\mu$.*

PROOF. (i) For bounded $\mu$, this is just the uniform version of (3). In general, we may write $\mu = \mu' + \mu''$ for bounded $\mu'$ and let $\xi = \xi' + \xi''$ be the corresponding decomposition of $\xi$. Then

$$P_\mu\{\xi_t B_x^\varepsilon > 0\} \leq P_\mu\{\xi_t' B_x^\varepsilon > 0\} + P_\mu\{\xi_t'' B_x^\varepsilon > 0\}$$
$$= P_{\mu'}\{\xi_t B_x^\varepsilon > 0\} + P_{\mu''}\{\xi_t B_x^\varepsilon > 0\},$$

and so by Lemmas 4.1 and 4.2(i)

$$|P_\mu\{\xi_t B_x^\varepsilon > 0\} - P_{\mu'}\{\xi_t B_x^\varepsilon > 0\}| \leq P_{\mu''}\{\xi_t B_x^\varepsilon > 0\}$$
$$\lesssim t\varepsilon^{d-2}(\mu'' * p_{t(\varepsilon)})(x).$$

For any $r > 0$ and for $\varepsilon_0 > 0$ small enough, there exists by Lemma 3.1 a $t' > 0$ such that
$$p_{t(\varepsilon)}(u-x) \lesssim p_{t'}(u), \qquad |x| \leq r, \ \varepsilon < \varepsilon_0, \ u \in \mathbb{R}^d,$$
which implies $(\mu'' * p_{t(\varepsilon)})(x) \leq \mu'' p_{t'}$ for the same $x$ and $\varepsilon$. Hence,
$$\|\varepsilon^{2-d} P_\mu\{\xi_t B^\varepsilon_\cdot > 0\} - c_d(\mu * p_t)\|_{B_0^r}$$
$$\lesssim \|\varepsilon^{2-d} P_{\mu'}\{\xi_t B^\varepsilon_\cdot > 0\} - c_d(\mu' * p_t)\| + \mu'' p_{t'},$$
which tends to 0 as $\varepsilon \to 0$ and then $\mu' \uparrow \mu$, by the result for bounded $\mu$ and dominated convergence.

(ii) First suppose that $\mu$ is bounded. Let $\varepsilon, h \to 0$ with $|\log h| \ll |\log \varepsilon|$, and write $\zeta_s$ for the ancestral process at time $s = t - h$. Then we get, uniformly on $\mathbb{R}^2$,
$$P_\mu\{\xi_t B^\varepsilon_\cdot > 0\} \approx E_\mu(\zeta_s * p^\varepsilon_h) = h^{-1} E_\mu(\xi_s * p^\varepsilon_h)$$
$$= h^{-1}(\mu * p_s * p^\varepsilon_h) \approx m(\varepsilon) |\log \varepsilon|^{-1} (\mu * p_s)$$
$$\approx m(\varepsilon) |\log \varepsilon|^{-1} (\mu * p_t).$$
To justify the first approximation, we see from Lemma 4.4(ii) that
$$|\log \varepsilon| \|P_\mu\{\xi_t B^\varepsilon_x > 0\} - E_\mu(\zeta_s * p^\varepsilon_h)\|$$
$$\lesssim \frac{|\log h| \|\mu * p_t\| + \|\mu * p_{t(h,\varepsilon)}\|^2}{|\log \varepsilon|} \lesssim \frac{|\log h|}{|\log \varepsilon|} \to 0.$$
For the second approximation, Lemma 5.2(ii) yields
$$\|h^{-1}|\log \varepsilon|(\mu * p_s * p^\varepsilon_h) - m(\varepsilon)(\mu * p_s)\|$$
$$\leq \|\mu\| \|h^{-1}|\log \varepsilon|(p_s * \tilde{p}^\varepsilon_h) - m(\varepsilon) p_s\| \to 0,$$
since the functions $p_s = p_{t-h}$ are uniformly bounded and equicontinuous for small $h > 0$. The third approximation holds since $m$ is bounded and
$$\|\mu * p_s - \mu * p_t\| \leq \|\mu\| \|p_s - p_t\| \to 0.$$
This completes the proof for bounded $\mu$. The extension to the general case may be accomplished by the same argument as for (i).

To prove the indicated version of (i) for the clusters $\eta_t$, we see from Lemmas 4.1 and 4.2(i) that
$$\varepsilon^{2-d} |t^{-1} P_\mu\{\eta_t B^\varepsilon_x > 0\} - P_\mu\{\xi_t B^\varepsilon_x > 0\}| \lesssim \varepsilon^{2-d}(t^{-1} P_\mu\{\eta_t B^\varepsilon_x > 0\})^2$$
$$\lesssim \varepsilon^{d-2}((\mu * p_{t(\varepsilon)})(x))^2.$$
For bounded $\mu$, this clearly tends to 0 as $\varepsilon \to 0$, uniformly in $x$. In general, Lemmas 3.1 and 3.2(iv) show that the right-hand side tends to 0, uniformly for bounded $x$. This proves the cluster version of (i), and the proof in case of (ii) is similar. $\square$

**6. Neighborhood measures.** For any measure $\mu$ on $\mathbb{R}^d$ and constant $\varepsilon > 0$, we define the associated *neighborhood measure* $\mu^\varepsilon$ as the restriction of Lebesgue measure $\lambda^d$ to the $\varepsilon$-neighborhood of $\operatorname{supp}\mu$, so that $\mu^\varepsilon$ has Lebesgue density $1\{\mu B_x^\varepsilon > 0\}$. In this section, we study the neighborhood measures of clusters $\eta_h$ associated with a DW-process in $\mathbb{R}^d$. This will prepare for the proof of the Lebesgue approximation of DW-processes in Section 7. We begin with some estimates of first and second moments.

LEMMA 6.1. *Let $\eta_1$ be the unit cluster of a DW-process in $\mathbb{R}^d$. Then as $\varepsilon \to 0$, we have:*

(i) $\|\varepsilon^{2-d} E_0 \eta_1^\varepsilon - c_d(p_1 \cdot \lambda^d)\| \to 0$ *for $d \geq 3$,*
(ii) $\||\log\varepsilon| E_0 \eta_1^\varepsilon - m(\varepsilon)(p_1 \cdot \lambda^2)\| \to 0$ *for $d = 2$,*
(iii) $E_0 \|\eta_1^\varepsilon\|^2 \asymp (E_0 \|\eta_1^\varepsilon\|)^2 \asymp \varepsilon^{2(d-2)}$ *for $d \geq 3$,*
(iv) $E_0 \|\eta_1^\varepsilon\|^2 \asymp (E_0 \|\eta_1^\varepsilon\|)^2 \asymp |\log\varepsilon|^{-2}$ *for $d = 2$.*

PROOF. (i) Fubini's theorem yields $E_0 \eta_1^\varepsilon = p_1^\varepsilon \cdot \lambda^d$, and so for $d \geq 3$

$$\|\varepsilon^{2-d} E_0 \eta_1^\varepsilon - c_d(p_1 \cdot \lambda^d)\| = \lambda^d |\varepsilon^{2-d} p_1^\varepsilon - c_d p_1|. \tag{6}$$

Here the integrand on the right tends to 0 as $\varepsilon \to 0$ by Theorem 5.3(i), and by Lemma 4.2(i) it is bounded by $C_d p_{1'} + c_d p_1 \to (C_d + c_d) p_1$ for some constant $C_d > 0$, where $1' = 1 + \varepsilon^2$. Since both sides have the same integral $C_d + c_d$, the integral in (6) tends to 0 by Theorem 1.21 in [13].

(ii) Use a similar argument based on Theorem 5.3(ii) and Lemma 4.2(ii).

(iii) For a DW-process $\xi$, let $\zeta_s$ be the process of ancestors of $\xi_1$ at time $s = 1 - h$, where $\varepsilon^2 \leq h \leq 1$, and denote the generated $h$-clusters by $\eta_h^i$. For any $x_1, x_2 \in \mathbb{R}^d$, write $x_i = \bar{x} \pm r$. Using Lemmas 3.3(i)–(ii) and 4.2(i), the conditional independence of the subclusters, the Cox property of $\zeta_s$ and the semigroup property of $p_t$, we obtain with $h' = h + \varepsilon^2$ and $1' = 1 + \varepsilon^2$

$$E_{\delta_0} \sum_{i \neq j} 1\{\eta_h^i B_{x_1}^\varepsilon \wedge \eta_h^j B_{x_2}^\varepsilon > 0\}$$

$$= \iint_{u_1 \neq u_2} p_h^\varepsilon(x_1 - u_1) p_h^\varepsilon(x_2 - u_2) E_{\delta_0} \zeta_s^2(du_1\, du_2)$$

$$\lesssim \varepsilon^{2(d-2)} \iint p_{h'}(x_1 - u_1) p_{h'}(x_2 - u_2) E_{\delta_0} \xi_s^2(du_1\, du_2)$$

$$= \varepsilon^{2(d-2)} ((p_s^{\otimes 2} + q_s) * p_{h'}^{\otimes 2})(x_1, x_2)$$

$$\leq \varepsilon^{2(d-2)} (p_{1'}^{\otimes 2} + q_{1'})(x_1, x_2)$$

$$\lesssim \varepsilon^{2(d-2)} p_{1'}(\bar{x}) p_{1'}(r) |r|^{2-d}.$$

Next we may combine the previously mentioned properties with Lemmas 3.2(iv) and 4.2(i), Cauchy's inequality, the parallelogram identity, and the special form of the densities $p_t$, to obtain

$$E_{\delta_0} \sum_i 1\{\eta_h^i B_{x_1}^\varepsilon \wedge \eta_h^i B_{x_2}^\varepsilon > 0\} = \int P_u\{\eta_h B_{x_1}^\varepsilon \wedge \eta_h B_{x_2}^\varepsilon > 0\} E_{\delta_0} \zeta_s(du)$$

$$\leq h^{-1} \int (p_h^\varepsilon(x_1-u) p_h^\varepsilon(x_2-u))^{1/2} E_{\delta_0} \xi_s(du)$$

$$\lesssim \varepsilon^{d-2} \int (p_{h'}(x_1-u) p_{h'}(x_2-u))^{1/2} p_s(u) \, du$$

$$\lesssim \varepsilon^{d-2} \int (p_{h'/2}(\bar{x}-u) p_{h'/2}(r))^{1/2} p_s(u) \, du$$

$$\lesssim \varepsilon^{d-2} h^{d/2} (p_{h'} * p_s)(\bar{x}) p_{h'}(r)$$

$$= \varepsilon^{d-2} h^{d/2} p_{1'}(\bar{x}) p_{h'}(r).$$

Since $\xi_1$ is the sum of $\kappa$ independent unit clusters, where $\kappa$ is Poisson under $P_{\delta_0}$ with mean 1, the previous estimates remain valid for the subclusters of $\eta$ of age $h$. Since $\eta_1^\varepsilon$ has Lebesgue density $1\{\eta_1 B_x^\varepsilon > 0\}$, Fubini's theorem yields

$$E_0 \|\eta_1^\varepsilon\|^2 = \iint P_0\{\eta_1 B_{x_i}^\varepsilon \wedge \eta_1 B_{x_i}^\varepsilon > 0\} \, dx_1 \, dx_2$$

$$\lesssim \iint dx_1 \, dx_2 \, E_{\delta_0} \sum_{i,j} 1\{\eta_h^i B_{x_1}^\varepsilon \wedge \eta_h^j B_{x_2}^\varepsilon > 0\}$$

$$\lesssim \iint (\varepsilon^{2(d-2)} p_{1'}(r) |r|^{2-d} + \varepsilon^{d-2} h^{d/2} p_{h'}(r)) p_{1'}(\bar{x}) \, d\bar{x} \, dr$$

$$\lesssim \varepsilon^{2(d-2)} + \varepsilon^{d-2} h^{d/2},$$

where, in the last step, we used the fact that

$$\int p_1(r) |r|^{2-d} \, dr \lesssim \int_0^\infty v e^{-v^2/2} \, dv < \infty.$$

Taking $h = \varepsilon^2$, we get by (i) and Jensen's inequality

$$\varepsilon^{2(d-2)} \asymp \|E_0 \eta_1^\varepsilon\|^2 \leq E_0 \|\eta_1^\varepsilon\|^2 \lesssim \varepsilon^{2(d-2)} + \varepsilon^{2d-2} \asymp \varepsilon^{2(d-2)}.$$

(iv) Suppose that $\varepsilon^2 \ll h \to 0$. Using Lemmas 3.3(iii) and 4.2(ii), we get as before

$$E_{\delta_0} \sum_{i \neq j} 1\{\eta_h^i B_{x_1}^\varepsilon \wedge \eta_h^j B_{x_2}^\varepsilon > 0\} \lesssim (\log(h/\varepsilon^2))^{-2} (p_{1'}^{\otimes 2} + q_{1'})(x_1, x_2)$$

$$\lesssim |\log \varepsilon|^{-2} p_{1'}(\bar{x}) p_1(r) \log(|r|^{-1} \vee e),$$

$$E_{\delta_0} \sum_i 1\{\eta_h^i B_{x_1}^\varepsilon \wedge \eta_h^i B_{x_2}^\varepsilon > 0\} \lesssim h |\log \varepsilon|^{-1} p_{1'}(\bar{x}) p_{h'}(r),$$

where $1' - 1 = h' - h \lesssim h|\log \varepsilon|^{-1/2}$. Noting that

$$\int p_1(r) \log(|r|^{-1} \vee e) \, dr \lesssim \int_{|r|e<1} |\log |r|| \, dr + \int p_1(r) \, dr < \infty,$$

we get by combination

$$E_0 \|\eta_1^\varepsilon\|^2 \lesssim \int\int (|\log \varepsilon|^{-2} p_1(r) \log(|r|^{-1} \vee e) + h|\log \varepsilon|^{-1} p_{h'}(r)) p_{1'}(\bar{x}) \, d\bar{x} \, dr$$

$$\lesssim |\log \varepsilon|^{-2} + h|\log \varepsilon|^{-1}.$$

Choosing $h = |\log \varepsilon|^{-1} \gg \varepsilon^2$ and combining with (ii) gives

$$|\log \varepsilon|^{-2} \asymp \|E_0 \eta_1^\varepsilon\|^2 \leq E_0 \|\eta_1^\varepsilon\|^2 \lesssim |\log \varepsilon|^{-2}. \qquad \square$$

This leads to some moment estimates for a Poisson "forest" of clusters. Recall that $p_h^\varepsilon(x) = P_x\{\eta_h B_0^\varepsilon > 0\}$ and write $(\eta_h^i)^\varepsilon = \eta_h^{i\varepsilon}$ for convenience.

LEMMA 6.2. *Let the $\eta_h^i$ be conditionally independent $h$-clusters in $\mathbb{R}^d$, rooted at the points of a Poisson process $\xi$ with $E\xi = \mu$. Fix any measurable function $f \geq 0$ on $\mathbb{R}^d$ and let $h \geq \varepsilon \to 0$. Then:*

(i) $E_\mu \sum_i \eta_h^{i\varepsilon} = (\mu * p_h^\varepsilon) \cdot \lambda^d$ *for $d \geq 2$,*
(ii) $\text{Var}_\mu \sum_i \eta_h^{i\varepsilon} f \lesssim h^2 \varepsilon^{2(d-2)} \|f\|^2 \|\mu\|$ *for $d \geq 3$,*
(iii) $\text{Var}_\mu \sum_i \eta_h^{i\varepsilon} f \lesssim h^2 |\log \varepsilon|^{-2} \|f\|^2 \|\mu\|$ *for $d = 2$.*

PROOF. (i) By Fubini's theorem and the definitions of $\eta_h^\varepsilon$ and $p_h^\varepsilon$, we have

$$E_x \eta_h^\varepsilon f = E_x \int 1\{\eta_h B_u^\varepsilon > 0\} f(u) \, du = (p_h^\varepsilon * f)(x),$$

and so by independence

$$(7) \qquad E\left[\sum_i \eta_h^{i\varepsilon} f \Big| \xi\right] = \int \xi(dx) E_x \eta_h^\varepsilon f = \xi(p_h^\varepsilon * f).$$

Hence, by Fubini's theorem

$$E_\mu \sum_i \eta_h^{i\varepsilon} f = E_\mu \xi(p_h^\varepsilon * f) = \mu(p_h^\varepsilon * f) = ((\mu * p_h^\varepsilon) \cdot \lambda^d) f.$$

(ii) By Lemma 2.1(ii) we have

$$\|\eta_h^\varepsilon\| = \int 1\{\eta_h B_x^\varepsilon > 0\} \, dx \stackrel{d}{=} \int 1\{\eta_1 B_{x/\sqrt{h}}^{\varepsilon/\sqrt{h}} > 0\} \, dx$$

$$= h^{d/2} \int 1\{\eta_1 B_x^{\varepsilon/\sqrt{h}} > 0\} \, dx = h^{d/2} \|\eta_1^{\varepsilon/\sqrt{h}}\|,$$

and so by Lemma 6.1(iii)

$$\text{Var}_x(\eta_h^\varepsilon f) \leq E_x(\eta_h^\varepsilon f)^2 \leq E\|\eta_h^\varepsilon\|^2 \|f\|^2 = h^d E\|\eta_1^{\varepsilon/\sqrt{h}}\|^2 \|f\|^2$$
$$\lesssim h^d (\varepsilon/\sqrt{h})^{2(d-2)} \|f\|^2 = \varepsilon^{2(d-2)} h^2 \|f\|^2.$$

Hence, by independence

$$E_\mu \text{Var}\left[\sum_i \eta_h^{i\varepsilon} f \Big| \xi\right] = E_\mu \int \xi(dx) \text{Var}_x(\eta_h^\varepsilon f)$$
$$\lesssim \varepsilon^{2(d-2)} h^2 \|f\|^2 \|\mu\|.$$

Since $\lambda^d p_h^\varepsilon \lesssim \varepsilon^{d-2} h$ by Lemma 4.2(i) and $\text{Var}_\mu(\xi f) = \mu f^2$, we get from (7)

$$\text{Var}_\mu E\left[\sum_i \eta_h^{i\varepsilon} f \Big| \xi\right] = \text{Var}_\mu \xi(p_h^\varepsilon * f) = \mu(p_h^\varepsilon * f)^2$$
$$\leq \|f\|^2 \|\mu\| (\lambda^d p_h^\varepsilon)^2 \lesssim \varepsilon^{2(d-2)} h^2 \|f\|^2 \|\mu\|.$$

Combining those estimates yields

$$\text{Var}_\mu \sum_i \eta_h^{i\varepsilon} f = E_\mu \text{Var}\left[\sum_i \eta_h^{i\varepsilon} f \Big| \xi\right] + \text{Var}_\mu E\left[\sum_i \eta_h^{i\varepsilon} f \Big| \xi\right]$$
$$\lesssim \varepsilon^{2(d-2)} h^2 \|f\|^2 \|\mu\|.$$

(iii) Since $h \geq \varepsilon$, we get by Lemma 6.1(iv)

$$\text{Var}_x(\eta_h^\varepsilon f) \lesssim h^2 |\log(\varepsilon/\sqrt{h})|^{-2} \|f\|^2 \lesssim h^2 |\log \varepsilon|^{-2} \|f\|^2,$$

and so

$$E_\mu \text{Var}\left[\sum_i \eta_h^{i\varepsilon} f \Big| \xi\right] \lesssim h^2 |\log \varepsilon|^{-2} \|f\|^2 \|\mu\|.$$

Next Lemma 4.2(ii) yields $\lambda^2 p_h^\varepsilon \lesssim h |\log \varepsilon|^{-1}$, so as before

$$\text{Var}_\mu E\left[\sum_i \eta_h^{i\varepsilon} f \Big| \xi\right] \lesssim h^2 |\log \varepsilon|^{-2} \|f\|^2 \|\mu\|.$$

The stated estimate now follows by combination. □

We also need to estimate the overlap between subclusters.

LEMMA 6.3. *Let $\xi$ be a DW-process in $\mathbb{R}^d$, and for fixed $t > 0$, let $\eta_h^i$ denote the subclusters in $\xi_t$ of age $h > 0$. Fix a $\mu \in \hat{\mathcal{M}}_d$. Then:*

(i) *for $d \geq 3$ and as $\varepsilon^2 \leq h \to 0$,*

$$E_\mu \left\| \sum_i \eta_h^{i\varepsilon} - \xi_t^\varepsilon \right\| \lesssim (\varepsilon^2/\sqrt{h})^{d-2},$$

(ii) *for $d = 2$ and as $\varepsilon \leq h \to 0$,*

$$E_\mu \left\| \sum_i \eta_h^{i\varepsilon} - \xi_t^\varepsilon \right\| \lesssim |\log h||\log \varepsilon|^{-2}.$$

PROOF. (i) Let $\kappa_h^\varepsilon(x)$ denote the number of subclusters of age $h$ hitting $B_x^\varepsilon$ at time $t$. Then Lemma 4.4(i) yields, with $t' = t + \varepsilon^2$,

$$\begin{aligned}
E_\mu \left\| \sum_i \eta_h^{i\varepsilon} - \xi_t^\varepsilon \right\| &= E_\mu \int \left| \sum_i 1\{\eta_h^i B_x^\varepsilon > 0\} - 1\{\xi B_x^\varepsilon > 0\} \right| dx \\
&= \int E_\mu (\kappa_h^\varepsilon(x) - 1)_+ \, dx \\
&\lesssim \varepsilon^{2(d-2)} \lambda^d (h^{1-d/2}(\mu * p_t) + (\mu * p_{t'})^2) \\
&\lesssim \varepsilon^{2(d-2)} (h^{1-d/2} \|\mu\| + t^{-d/2} \|\mu\|^2).
\end{aligned}$$

(ii) Using Lemma 4.4(ii), we get instead

$$\begin{aligned}
E_\mu \left\| \sum_i \eta_h^{i\varepsilon} - \xi_t^\varepsilon \right\| &\lesssim |\log \varepsilon|^{-2} \lambda^2 (\log(t/h)(\mu * p_t) + (\mu * p_{t'})^2) \\
&\lesssim |\log \varepsilon|^{-2} (|\log h| \|\mu\| + t^{-1} \|\mu\|^2),
\end{aligned}$$

for a suitable choice of $t' \geq t$. □

**7. Lebesgue approximation.** Given a DW-process $\xi$ in $\mathbb{R}^d$, we prove for any $d \geq 2$ and for fixed $t > 0$ that $\xi_t$ can be approximated, both a.s. and in $L^1$, by suitably normalized versions of the neighborhood measures $\xi_t^\varepsilon$, as defined in Section 6. For $d \geq 3$, this result is essentially due to Tribe [27]. Write $\tilde{c}_d = 1/c_d$ and $\tilde{m} = 1/m$ for convenience, where $c_d$ and $m$ are such as in Section 5.

THEOREM 7.1. *Let the DW-process $\xi$ in $\mathbb{R}^d$ be locally finite under $P_\mu$, and fix a $t > 0$. Then under $P_\mu$, we have as $\varepsilon \to 0$:*

(i) $\tilde{c}_d \varepsilon^{2-d} \xi_t^\varepsilon \xrightarrow{v} \xi_t$ *a.s. and in $L^1$ for $d \geq 3$,*
(ii) $\tilde{m}(\varepsilon)|\log \varepsilon| \xi_t^\varepsilon \xrightarrow{v} \xi_t$ *a.s. and in $L^1$ for $d = 2$.*

*This remains true in the weak sense when $\mu$ is bounded. The weak versions hold even for the clusters $\eta_t$ when $\|\mu\| = 1$.*

PROOF. We use a new approach, explained in detail only for $d \geq 3$. (i) Let $d \geq 3$, and fix any $t > 0$, $\mu \in \hat{\mathcal{M}}_d$ and $f \in C_K^d$. Write $\eta_h^i$ for the subclusters of $\xi_t$ of age $h$. Since the ancestors of $\xi_t$ at time $s = t - h$ form a Cox process directed by $\xi_s/h$, Lemma 6.2(i) yields

$$E_\mu\left[\sum_i \eta_h^{i\varepsilon} f \,\Big|\, \xi_s\right] = h^{-1}\xi_s(p_h^\varepsilon * f),$$

and so by Lemma 6.2(ii)

$$E_\mu\left|\sum_i \eta_h^{i\varepsilon} f - h^{-1}\xi_s(p_h^\varepsilon * f)\right|^2 = E_\mu \operatorname{Var}\left[\sum_i \eta_h^{i\varepsilon} f \,\Big|\, \xi_s\right]$$
$$\lesssim \varepsilon^{2(d-2)} h^2 \|f\|^2 E_\mu \|\xi_s/h\|$$
$$= \varepsilon^{2(d-2)} h \|f\|^2 \|\mu\|.$$

Combining with Lemma 6.3(i) gives

$$E_\mu |\xi_t^\varepsilon f - h^{-1}\xi_s(p_h^\varepsilon * f)|$$
$$\leq E_\mu\left|\xi_t^\varepsilon f - \sum_i \eta_h^{i\varepsilon} f\right| + E_\mu\left|\sum_i \eta_h^{i\varepsilon} f - h^{-1}\xi_s(p_h^\varepsilon * f)\right|$$
$$\lesssim \varepsilon^{2(d-2)} h^{1-d/2} \|f\| + \varepsilon^{d-2} h^{1/2} \|f\|$$
$$= \varepsilon^{d-2}(\sqrt{h} + (\varepsilon/\sqrt{h})^{d-2})\|f\|.$$

Taking $h = \varepsilon = r^n$ for a fixed $r \in (0,1)$ and writing $s_n = t - r^n$, we obtain

$$E_\mu \sum_n r^{n(2-d)} |\xi_t^{r^n} f - r^{-n}\xi_{s_n}(p_{r^n}^{r^n} * f)| \lesssim \sum_n (r^{n/2} + r^{n(d-2)/2})\|f\| < \infty,$$

which implies

$$(8) \qquad r^{n(2-d)} |\xi_t^{r^n} f - r^{-n}\xi_{s_n}(p_{r^n}^{r^n} * f)| \to 0 \qquad \text{a.s. } P_\mu.$$

Now we write

$$|\varepsilon^{2-d}\xi_t^\varepsilon f - c_d \xi_t f| \leq \varepsilon^{2-d}|\xi_t^\varepsilon f - h^{-1}\xi_s(p_h^\varepsilon * f)| + c_d |\xi_s f - \xi_t f|$$
$$+ \|\xi_s\|\|\varepsilon^{2-d} h^{-1}(p_h^\varepsilon * f) - c_d f\|.$$

Using (8), Lemma 5.2(i) and the a.s. weak continuity of $\xi$ (cf. Proposition 2.15 in [8]), we see that the right-hand side tends a.s. to 0 as $n \to \infty$, which implies $\varepsilon^{2-d}\xi_t^\varepsilon f \to c_d \xi_t f$ a.s. as $\varepsilon \to 0$ along the sequence $(r^n)$ for any fixed $r \in (0,1)$. Since this holds simultaneously, outside a fixed null set, for all rational $r \in (0,1)$, the a.s. convergence extends by Lemma 2.3 to the entire interval $(0,1)$.

Now let $\mu \in \mathcal{M}_d$ be arbitrary with $\mu p_t < \infty$ for all $t > 0$. Write $\mu = \mu' + \mu''$ for bounded $\mu'$, and let $\xi = \xi' + \xi''$ be the corresponding decomposition of $\xi$ into independent components with initial measures $\mu'$ and $\mu''$. Fixing an $r > 1$ with $\operatorname{supp} f \subset B_0^{r-1}$ and using the result for bounded $\mu$, we get a.s. on $\{\xi_t'' B_0^r = 0\}$

$$\varepsilon^{2-d}\xi_t^\varepsilon f = \varepsilon^{2-d}\xi_t'^\varepsilon f \to c_d \xi_t' f = c_d \xi_t f.$$

As $\mu' \uparrow \mu$, we get by Lemma 4.3

$$P_\mu\{\xi_t'' B_0^r = 0\} = P_{\mu''}\{\xi_t B_0^r = 0\} \to 1,$$

and the a.s. convergence extends to $\mu$. Applying this result to a countable, convergence-determining class of functions $f$ (cf. Lemma 3.2.1 in [3]), we obtain the required a.s. vague convergence. If $\mu$ is bounded, then $\xi_t$ has a.s. bounded support (cf. Corollary 6.8 in [8]), and the a.s. convergence remains valid in the weak sense.

To prove the convergence in $L^1$, we note that for any $f \in C_K^d$

$$\varepsilon^{2-d} E_\mu \xi_t^\varepsilon f = \varepsilon^{2-d} \int P_\mu\{\xi_t B_x^\varepsilon > 0\} f(x)\, dx$$

(9)

$$\to \int c_d (\mu * p_t)(x) f(x)\, dx = c_d E_\mu \xi_t f,$$

by Theorem 5.3(i). Combining this with the a.s. convergence under $P_\mu$ and using Proposition 4.12 in [13], we obtain $E_\mu|\varepsilon^{2-d}\xi_t^\varepsilon f - c_d \xi_t f| \to 0$. For bounded $\mu$, (9) extends to any $f \in C_b^d$ by dominated convergence based on Lemmas 4.1 and 4.2(i), together with the fact that $\lambda^d(\mu * p_t) = \|\mu\| < \infty$ by Fubini's theorem.

(ii) Let $d = 2$, and fix any $t$, $\mu$ and $f$ as before. Using Lemma 6.2(iii), we see as in part (i) that

$$E_\mu \left|\sum_i \eta_h^{i\varepsilon} f - h^{-1}\xi_s(p_h^\varepsilon * f)\right|^2 \lesssim h|\log \varepsilon|^{-2} \|f\|^2 \|\mu\|.$$

Combining with Lemma 6.3(ii), we now get for fixed $\mu$ and $f$

$$E_\mu |\xi_t^\varepsilon f - h^{-1}\xi_s(p_h^\varepsilon * f)| \lesssim h^{1/2}|\log\varepsilon|^{-1} + |\log h||\log\varepsilon|^{-2}.$$

Choosing $\sqrt{h} = |\log\varepsilon|^{-1} = r^n$ for a fixed $r \in (0,1)$, we get

(10)
$$|\log\varepsilon| E_\mu |\xi_t^\varepsilon f - h^{-1}\xi_s(p_h^\varepsilon * f)| \lesssim h^{1/2} + |\log h||\log\varepsilon|^{-1}$$
$$= r^n + 2n|\log r|r^n \lesssim r^{n/2}.$$

Now we write

$$|\tilde{m}(\varepsilon)|\log\varepsilon|\xi_t^\varepsilon f - \xi_t f| \leq |\log\varepsilon||\xi_t^\varepsilon f - h^{-1}\xi_s(p_h^\varepsilon * f)| + |\xi_s f - \xi_t f|$$
$$+ \|\xi_s\|\|h^{-1}\tilde{m}(\varepsilon)|\log\varepsilon|(p_h^\varepsilon * f) - f\|.$$

Letting $\sqrt{h} = |\log \varepsilon|^{-1} = r^n$ with $n \to \infty$, we see from (10), Lemma 5.2(ii) and the weak continuity of $\xi$ that the right-hand side tends to 0 a.s. Writing $\varepsilon = e^{-1/s}$ and putting $\tilde{\xi}_t^s = \xi_t^\varepsilon$, we conclude that

$$\tilde{m}(e^{-1/s})s^{-1}\tilde{\xi}_t^s f \to \xi_t f \qquad \text{a.s. } P_\mu, \tag{11}$$

as $s \to 0$ along $(r^n)$ for any $r \in (0,1)$. Since the function $t \mapsto \log \tilde{m}(\exp(-e^t))$ is bounded and uniformly continuous on $\mathbb{R}_+$ by Lemma 5.1, (11) remains true along $(0,1)$ by Lemma 2.3. Hence, $\tilde{m}(\varepsilon)|\log \varepsilon|\xi_t^\varepsilon f \to \xi_t f$ a.s. $P_\mu$ for fixed $f$ and bounded $\mu$, which extends as before to $\tilde{m}(\varepsilon)|\log \varepsilon|\xi_t^\varepsilon \xrightarrow{v} \xi_t$ a.s., even when $\mu$ is unbounded.

To prove the corresponding $L^1$-convergence, let $f \in C_K^d$ and conclude from Theorem 5.3(ii) that

$$\tilde{m}(\varepsilon)|\log \varepsilon|E_\mu \xi_t^\varepsilon f = \tilde{m}(\varepsilon)|\log \varepsilon|\int P_\mu\{\xi_t B_x^\varepsilon > 0\}f(x)\,dx$$

$$\to \int (\mu * p_t)(x)f(x)\,dx = E_\mu \xi_t f.$$

For bounded $\mu$, this extends by dominated convergence to any $f \in C_b^d$. The assertion now follows as before by combination with the corresponding a.s. convergence.

To extend (i) and (ii) to the individual clusters $\eta_t$, let $\zeta_0$ denote the process of ancestors of $\xi_t$ at time 0, and note that

$$P_0\{\eta_t \in \cdot\} = P_{\delta_0}[\xi_t \in \cdot |\|\zeta_0\| = 1],$$

where $P_{\delta_0}\{\|\zeta_0\| = 1\} = t^{-1}e^{-1/t} > 0$. The a.s. convergence then follows from the corresponding statement for $\xi_t$. To obtain the weak $L^1$-convergence in this case, we note that for $f \in C_b^d$ and $d \geq 3$ or $d = 2$, respectively,

$$\varepsilon^{2-d}E_0\eta_t^\varepsilon f = \varepsilon^{2-d}\lambda^d(p_t^\varepsilon f) \to c_d t \lambda^d(p_t f) = c_d E_0 \eta_t f,$$

$$\tilde{m}(\varepsilon)|\log \varepsilon|E_0\eta_t^\varepsilon f = \tilde{m}(\varepsilon)|\log \varepsilon|\lambda^d(p_t^\varepsilon f) \to t\lambda^d(p_t f) = E_0 \eta_t f,$$

by dominated convergence based on Lemma 4.2 and Theorem 5.3. $\square$

For the intensity measures in Theorem 7.1, we have even convergence in total variation.

COROLLARY 7.2. *Let $\xi$ be a DW-process in $\mathbb{R}^d$. Then for any $t > 0$ and bounded $\mu$, we have as $\varepsilon \to 0$:*

(i) $\|\varepsilon^{2-d}E_\mu\xi_t^\varepsilon - c_d E_\mu \xi_t\| \to 0$ *for $d \geq 3$,*
(ii) $\||\log \varepsilon|E_\mu\xi_t^\varepsilon - m(\varepsilon)E_\mu \xi_t\| \to 0$ *for $d = 2$.*

*The results remain true for the clusters $\eta_t$, and they also hold locally for $\xi_t$ whenever $\xi$ is locally finite under $P_\mu$.*

PROOF. The two conditions are equivalent to the statements

$$\int |\varepsilon^{2-d} P_\mu\{\xi_t B_x^\varepsilon > 0\} - c_d(\mu * p_t)(x)|\, dx \to 0,$$

$$\int ||\log \varepsilon| P_\mu\{\xi_t B_x^\varepsilon > 0\} - m(\varepsilon)(\mu * p_t)|\, dx \to 0,$$

which are $L^1$-versions of Theorem 5.3 and follow as before by dominated convergence. $\square$

**8. Strong approximation for $d \geq 3$.** Here we prove that the distribution of a DW-process of dimension $d \geq 3$ admits a local approximation, in the sense of total variation, by a stationary and self-similar pseudo-random measure $\tilde{\xi}$. A related but weaker result is mentioned without proof in [5], page 119, with reference to some unpublished work with Iscoe.

For any $B \in \hat{\mathcal{B}}^d$ we write $\|\cdot\|_B$ for the total variation on the set $H_B = \{\mu; \mu B > 0\}$, equipped with the $\sigma$-field $\mathcal{H}_B$ generated by the restriction map $\mu \mapsto 1_B \cdot \mu$.

THEOREM 8.1. *For $d \geq 3$, let the DW-process $\xi$ in $\mathbb{R}^d$ be locally finite under $P_\mu$. Then there exists a pseudo-random measure $\tilde{\xi}$ on $\mathbb{R}^d$ such that:*

(i) *as $\varepsilon \to 0$ for fixed $B \in \hat{\mathcal{B}}^d$ and $t > 0$,*

$$\|\varepsilon^{2-d} P_\mu\{\varepsilon^{-2}\xi_t S_\varepsilon \in \cdot\} - \mu p_t \tilde{P}\{\tilde{\xi} \in \cdot\}\|_B \to 0,$$

*and similarly for the clusters $\eta_t$ with $p_t$ replaced by $tp_t$,*

(ii) *for any $r > 0$ and $a \in \mathbb{R}^d$,*

$$\tilde{P}\{\tilde{\xi} S_r \theta_a \in \cdot\} = r^{d-2} \tilde{P}\{r^2 \tilde{\xi} \in \cdot\}.$$

PROOF. Fix any $t > h > 0$ and $B \in \hat{\mathcal{B}}^d$, and consider any $\mathcal{H}_B$-measurable function $f \geq 0$ on $\mathcal{M}_d$ with $f \leq 1_{H_B}$. Consider the process $\zeta_s$ of ancestors of $\xi_t$ at time $s = t - h$, and let $\eta_h^i$ denote the associated $h$-clusters. As $h \to 0$ and $r = \varepsilon/\sqrt{h} \to 0$, we have the following chain of relations, explained in further detail below:

$$
\begin{aligned}
E_\mu f(\varepsilon^{-2}\xi_t S_\varepsilon) &= E_\mu f\left(\varepsilon^{-2} \sum_i \eta_h^i S_\varepsilon\right) \approx E_\mu \sum_i f(\varepsilon^{-2}\eta_h^i S_\varepsilon) \\
&= \int E_x f(\varepsilon^{-2}\eta_h S_\varepsilon) E_\mu \zeta_s(dx) \\
&= h^{-1} \int \mu(dy) \int p_s(x-y) E_x f(\varepsilon^{-2}\eta_h S_\varepsilon)\, dx \\
&\approx h^{-1} \mu p_t \int E_x f(\varepsilon^{-2}\eta_h S_\varepsilon)\, dx
\end{aligned}
$$

(12)

$$= h^{-1}\mu p_t \int E_{x/\sqrt{h}} f(h\varepsilon^{-2}\eta_1 S_{\varepsilon/\sqrt{h}})\,dx$$

$$= (\varepsilon/r)^{d-2}\mu p_t \int E_x f(r^{-2}\eta_1 S_r)\,dx.$$

Here the third relation holds by the conditional independence of the clusters, the fourth relation holds since $E_\mu \zeta_s = h^{-1} E_\mu \xi_s = h^{-1}(\mu * p_s) \cdot \lambda^d$, and the sixth relation holds by Lemma 2.1.

To justify the first approximation in (12), define $\kappa_h^\varepsilon$ as in Lemma 4.4 and fix a $c > 0$ with $B \subset B_0^c$. Then the mentioned lemma yields

$$\varepsilon^{2-d} E_\mu \left| f\left(\varepsilon^{-2}\sum_i \eta_h^i S_\varepsilon\right) - \sum_i f(\varepsilon^{-2}\eta_h^i S_\varepsilon) \right|$$

(13)
$$\leq \varepsilon^{2-d} E_\mu[\kappa_h^{c\varepsilon}; \kappa_h^{c\varepsilon} > 1]$$
$$\lesssim \varepsilon^{2-d}(c\varepsilon)^{2(d-2)}(h^{1-d/2}\mu p_t + (\mu p_{t(c\varepsilon)})^2) \lesssim r^{d-2} \to 0.$$

The second approximation in (12) amounts to replacing $p_s(x-y)$ by $p_t(y)$ in the inner integral. To estimate the resulting error, we note that by Lemma 4.2(i)

$$\varepsilon^{2-d} h^{-1} \left| \int \mu(dy) \int (p_s(y-x) - p_t(y)) E_x f(\varepsilon^{-2}\eta_h S_\varepsilon)\,dx \right|$$

(14)
$$\lesssim \int \mu(dy) \int |p_s(y-x) - p_t(y)| p_{h(c\varepsilon)}(x)\,dx$$
$$= \int \mu(dy) E|p_s(y - \gamma h_{c\varepsilon}^{1/2}) - p_t(y)|,$$

where $h_\varepsilon = h + \varepsilon^2$ and $\gamma$ denotes a standard normal random vector in $\mathbb{R}^d$. As $\varepsilon^2 \leq h \to 0$, we get $p_s(y - \gamma h_\varepsilon^{1/2}) \to p_t(y)$ a.s. by the joint continuity of $(x,t) \mapsto p_t(x)$. Since also

$$E p_s(y - \gamma h_{c\varepsilon}^{1/2}) = (p_s * p_{h(c\varepsilon)})(y) = p_{t(c\varepsilon)}(y) \to p_t(y),$$

the last expectation in (14) tends to 0 by Lemma 1.32 in [13]. Finally, since

$$E|p_s(y - \gamma h_{c\varepsilon}^{1/2}) - p_t(y)| \leq p_{t(c\varepsilon)}(y) + p_t(y) \lesssim p_{2t}(y),$$

where $\mu p_{2t} < \infty$, the right-hand side of (14) tends to 0 by dominated convergence.

This proves that, as $\varepsilon \ll r \to 0$ for fixed $\mu \in \mathcal{M}_d$, $B \in \hat{\mathcal{B}}^d$ and $t > 0$,

(15) $$\left\| \varepsilon^{2-d} P_\mu\{\varepsilon^{-2}\xi_t S_\varepsilon \in \cdot\} - r^{2-d}\mu p_t \int P_x\{r^{-2}\eta_1 S_r \in \cdot\}\,dx \right\|_B \to 0.$$

In particular, the first term on the left is uniformly Cauchy convergent on $H_B$ as $\varepsilon \to 0$. Hence, both terms converge as $\varepsilon \to 0$ and $r \to 0$, respectively, to a common limit of the form $\mu p_t \varphi_B$, where the set function $\varphi_B$ on $H_B$ is independent of $\mu$ and $t$. Thus,

$$\|\varepsilon^{2-d} P_\mu\{\varepsilon^{-2} \xi_t S_\varepsilon \in \cdot\} - \mu p_t \varphi_B\|_B \to 0, \tag{16}$$

where the uniformity of the convergence ensures that $\varphi_B$ is a bounded measure on $(H_B, \mathcal{H}_B)$.

Comparing the statements (16) for different sets $B$, we see that the $\varphi_B$ are all restrictions of a common set function $\varphi$ on $\bigcup_B \mathcal{H}_B$. We need to prove that $\varphi$ can be extended to a measure $\hat\varphi$ on $\bigcup_B H_B = \{\mu \in \mathcal{M}_d; \mu \neq 0\} = \mathcal{M}'_d$, endowed with the $\sigma$-field $\mathcal{H} = \bigvee_B \mathcal{H}_B$ generated by all projection maps $\mu \mapsto \mu B$. Choosing $\tilde P = \hat\varphi$ and letting $\tilde\xi$ denote the identity map on $\mathcal{M}'_d$, we may then write (16) in the form (i).

To construct $\hat\varphi$, it is enough for every fixed $B \in \hat{\mathcal{B}}^d$ to form the restriction $\hat\varphi_B$ of $\hat\varphi$ to $H_B$ with the trace $\sigma$-field $H_B \cap \mathcal{H}$, since the measure $\hat\varphi = \sup_B \hat\varphi_B$ has then the required properties. Writing $S = \mathcal{M}_d$ and $S_n = \mathcal{M}(B_0^n)$, for all $n$ satisfying $B_0^n \supset B$, we introduce the restriction maps $\pi_n \colon S \to S_n$ and $\pi_{n,k} \colon S_n \to S_k$, $n \geq k$. Put $\varphi'_n = \varphi_{B_0^n}(H_B \cap \cdot)$ and form the bounded measures $\psi_n = \varphi'_n \circ \pi_n^{-1}$ on $S_n$. Since $\psi_n \circ \pi_{n,k}^{-1} = \psi_k$ for all $n \geq k$, and since measures in $\mathcal{M}_d$ are measurably determined by their restrictions to the balls $B_0^n$, there exists by Corollary 6.15 in [13] a measure $\psi$ on $S$ with $\psi_n = \psi \circ \pi_n^{-1}$ for all $n$. Since the $\psi_n$ are restricted to $H_B$, so is $\psi$, and we see that $\hat\varphi_B = \psi$ has the desired properties.

To show that (i) remains true for the clusters $\eta_t$ with $p_t$ replaced by $tp_t$, we may apply the first four relations in (12)—as justified by (13)—with $h = t$ and $s = 0$, to get as $\varepsilon \to 0$ for fixed $B \in \hat{\mathcal{B}}^d$

$$\|t\varepsilon^{2-d} P_\mu\{\varepsilon^{-2} \xi_1 S_\varepsilon \in \cdot\} - \varepsilon^{2-d} P_\mu\{\varepsilon^{-2} \eta_t S_\varepsilon \in \cdot\}\|_B \to 0.$$

The required convergence now follows from (i).

To prove (ii), we may use the shift and semigroup properties of the operators $S_x$ and the shift invariance of $\lambda^d$ to get, for any $r, \varepsilon > 0$ and $a \in \mathbb{R}^d$,

$$\varepsilon^{2-d} \int P_x\{\varepsilon^{-2} \eta_1 S_\varepsilon S_r \theta_a \in \cdot\}\,dx = r^{d-2}(r\varepsilon)^{2-d} \int P_x\{r^2(r\varepsilon)^{-2} \eta_1 S_{r\varepsilon} \in \cdot\}\,dx.$$

Letting $\varepsilon \to 0$ for fixed $r$ and applying the cluster version of (i) to each side, we obtain (ii) on $(H_B, \mathcal{H}_B)$ for every $B \in \hat{\mathcal{B}}^d$, and the general result follows by a monotone class argument. □

The previous convergence extends to the associated Palm distributions, which will be useful in the next section.

THEOREM 8.2. *For $d \geq 3$, let the DW-process $\xi$ in $\mathbb{R}^d$ be locally finite under $P_\mu$, and let $\tilde\xi$ be such as in Theorem 8.1. Then $\tilde E \tilde\xi = \lambda^d$, and we may* *introduce the associated Palm distributions $P_\mu^0$ and $\tilde{P}^0$. Letting $\varepsilon \to 0$ for fixed $B \in \hat{\mathcal{B}}^d$ and $t > 0$, we have*

$$\|P_\mu^0\{\varepsilon^{-2}\xi_t S_\varepsilon \in \cdot\} - \tilde{P}^0\{\tilde{\xi} \in \cdot\}\|_B \to 0,$$

*and similarly with $\xi_t$ replaced by $\eta_t$.*

PROOF. Noting that $E_\mu \xi_t = t^{-1} E_\mu \eta_t = (\mu * p_t) \cdot \lambda^d$ and using the continuity in Lemma 3.2(iv), we get as $\varepsilon \to 0$ for fixed $B \in \hat{\mathcal{B}}^d$

(17) $$\varepsilon^{-d} E_\mu \xi_t(\varepsilon B) = t^{-1}\varepsilon^{-d} E_\mu \eta_t(\varepsilon B) \to \mu p_t \lambda^d B.$$

Using Lemma 3.3(i) above and Hint (2) in [26], page 239, we obtain

$$\mathrm{Var}_\mu \xi_t B_0^\varepsilon \lesssim E_\mu \xi_t B_0^\varepsilon \int_0^t (\varepsilon^d s^{-d/2} \wedge 1)\, ds$$
$$\lesssim \varepsilon^d \mu p_t \lambda^d B_0^1 \left( \int_0^{\varepsilon^2} ds + \varepsilon^d \int_{\varepsilon^2}^t s^{-d/2}\, ds \right) \lesssim \varepsilon^{d+2} \mu p_t.$$

Combining with (17) and Theorem 5.3(i), we get

(18) $$E_\mu[(\varepsilon^{-2}\xi_t B_0^\varepsilon)^2 | \xi_t B_0^\varepsilon > 0] = \frac{(\varepsilon^{-2} E_\mu \xi_t B_0^\varepsilon)^2 + \mathrm{Var}_\mu(\varepsilon^{-2}\xi_t B_0^\varepsilon)}{P_\mu\{\xi_t B_0^\varepsilon > 0\}}$$
$$\lesssim \frac{(\varepsilon^{d-2}\mu p_t)^2 + \varepsilon^{d-2}\mu p_t}{\varepsilon^{d-2}\mu p_t} \lesssim 1.$$

Next we see from Theorem 8.1(i) that, for $B_0^1 \subset B \in \hat{\mathcal{B}}^d$,

(19) $$\|P_\mu[\varepsilon^{-2}\xi_t S_\varepsilon \in \cdot | \xi_t B_0^\varepsilon > 0] - \tilde{P}[\tilde{\xi} \in \cdot | \tilde{\xi} B_0^1 > 0]\|_B \to 0.$$

By (18) the random variables $\varepsilon^{-2}\xi_t B_0^\varepsilon$ are uniformly integrable, conditionally on $\xi_t B_0^\varepsilon > 0$. Hence, by a uniform version of Lemma 4.11 in [13], we may extend (19) to

$$\|E_\mu[\varepsilon^{-2}\xi_t B_0^\varepsilon; \varepsilon^{-2}\xi_t S_\varepsilon \in \cdot | \xi_t B_0^\varepsilon > 0] - \tilde{E}[\tilde{\xi} B_0^1; \tilde{\xi} \in \cdot | \tilde{\xi} B_0^1 > 0]\|_B \to 0.$$

Combining this with Theorem 8.1(i) yields

(20) $$\|\varepsilon^{2-d} E_\mu[\varepsilon^{-2}\xi_t B_0^\varepsilon; \varepsilon^{-2}\xi_t S_\varepsilon \in \cdot] - \mu p_t \tilde{E}[\tilde{\xi} B_0^1; \tilde{\xi} \in \cdot]\|_B \to 0.$$

Since $B_0^1 \subset B$, we see from (17) and (20) that

(21) $$t^{-1}\varepsilon^{-d} E_\mu \eta_t B_0^\varepsilon = \varepsilon^{-d} E_\mu \xi_t B_0^\varepsilon \to \mu p_t \tilde{E}\tilde{\xi}B_0^1 = \mu p_t \lambda^d B_0^1.$$

Hence, by stationarity $\tilde{E}\tilde{\xi} = \lambda^d$, which justifies the definition of $\tilde{P}^0$. From (21) and Theorem 8.1(i) we obtain

$$E_\mu[\varepsilon^{-2}\eta_t B_0^\varepsilon | \eta_t B_0^\varepsilon > 0] \to \tilde{E}[\tilde{\xi}B_0^1 | \tilde{\xi} B_0^1 > 0],$$
$$\|P_\mu[\varepsilon^{-2}\eta_t S_\varepsilon \in \cdot | \eta_t B_0^\varepsilon > 0] - \tilde{P}[\tilde{\xi} \in \cdot | \tilde{\xi} B_0^1 > 0]\|_B \to 0.$$

By Lemma 4.11 of [13], now used in the opposite direction, we conclude that the variables $\varepsilon^{-2}\eta_t B_0^\varepsilon$ are uniformly integrable, conditionally on $\eta_t B_0^\varepsilon > 0$. Hence, by the uniform version of the same lemma

$$(22) \qquad \|\varepsilon^{2-d} E_\mu[\varepsilon^{-2}\eta_t B_0^\varepsilon; \varepsilon^{-2}\eta_t S_\varepsilon \in \cdot] - t\mu p_t \tilde{E}[\tilde{\xi} B_0^1; \tilde{\xi} \in \cdot]\|_B \to 0.$$

The asserted convergence now follows by Lemma 2.2, adapted to the case of a pseudo-random limiting measure $\tilde{\xi}$ with $\tilde{P}\{\tilde{\xi}B > 0\} < \infty$. Here conditions (i) and (ii) hold by (20) and (22), and Lemma 3.5 yields (iii) for the shifted Palm distributions of $\xi_t$ and $\eta_t$, based on an arbitrary initial measure $\mu$. □

**9. Local invariance for $d = 2$.** In two dimensions, the DW-process exhibits a completely different local behavior. Here we show that the measures $\xi_t$ at fixed times $t > 0$ are then locally invariant in a number of different ways. It is interesting to compare with the *diffusive clustering* discussed by Klenke [19].

THEOREM 9.1. *Let the DW-process $\xi$ in $\mathbb{R}^2$ be locally finite under $P_\mu$, and define $\rho_t^\varepsilon = \xi_t B_0^\varepsilon / \pi$ and $P_\mu^\varepsilon = P_\mu[\cdot|\rho_t^\varepsilon > 0]$. Then as $\varepsilon \to 0$ for fixed $t > 0$, we have:*

(i) $\xi_t S_\varepsilon / \rho_t^\varepsilon \xrightarrow{d} \lambda^2$ *under $P_\mu^0$,*
(ii) $E_\mu^\varepsilon |\xi_t S_\varepsilon f - \rho_t^\varepsilon \lambda^2 f| / E_\mu^\varepsilon \rho_t^\varepsilon \to 0$ *for all $f \in C_K^2$,*
(iii) $\mathrm{supp}(\xi_t S_\varepsilon) \xrightarrow{d} \mathbb{R}^2$ *under $P_\mu^\varepsilon$.*

*All statements remain true for the clusters $\eta_t$ when $\|\mu\| = 1$.*

Here (iii) means that $P_\mu^\varepsilon\{\xi_t S_\varepsilon B > 0\} \to 1$ for all open sets $B$. From (i) we see that $\mathrm{supp}(\xi_t S_\varepsilon) \xrightarrow{d} \mathbb{R}^2$ holds even under $P_\mu^0$. Statement (ii) is substantial, since the variables $\rho_t^\varepsilon$ are uniformly integrable under $P_\mu^\varepsilon$. However, it is not strong enough to imply (iii), and it is not clear whether (ii) can be strengthened to $\xi_t S_\varepsilon / \rho_t^\varepsilon \xrightarrow{d} \lambda^2$ under $P_\mu^\varepsilon$.

PROOF. We consider only $\xi_t$, the proof for $\eta_t$ being similar.
(i) Here Lemma 3.4 yields

$$(23) \qquad E_\mu(\xi_t S_\varepsilon f - \rho_t^\varepsilon \lambda^2 f)^2 \ll \varepsilon^4 |\log \varepsilon| \mu p_t, \qquad f \in C_K^2.$$

Using Cauchy's inequality and Lemma 4.2(ii), we get for fixed $f \in C_K^2$ and $r \in (0,1)$

$$(24) \qquad \begin{aligned} E_\mu \rho_t^{r\varepsilon} |\xi_t S_\varepsilon f / \rho_t^\varepsilon - \lambda^2 f| &\leq (E_\mu(\xi_t S_\varepsilon f - \rho_t^\varepsilon \lambda^2 f)^2 P_\mu\{\rho_t^\varepsilon > 0\})^{1/2} \\ &\ll (\varepsilon^4 |\log \varepsilon| \mu p_t |\log \varepsilon|^{-1} \mu p_t)^{1/2} \\ &= \varepsilon^2 \mu p_t \lesssim E_\mu \rho_t^{r\varepsilon}. \end{aligned}$$

Now define
$$f_r^+(x) = \sup_{|u|\leq r} f(x+u), \qquad f_r^-(x) = \inf_{|x|\leq r} f(x+u), \qquad x \in \mathbb{R}^2, \ r > 0,$$

and note that

$$E_\mu \rho_t^{r\varepsilon} \inf_{|x|\leq r\varepsilon} E_{\mu\theta_x}^0(|\xi_t S_\varepsilon f/\rho_t^\varepsilon - \lambda^2 f| \wedge 1)$$

$$\leq \int_{|x|\leq r\varepsilon} E_\mu \xi_t(dx) E_{\mu\theta_{-x}}^0 |\xi_t S_\varepsilon f/\rho_t^\varepsilon - \lambda^2 f|$$

$$\leq E_\mu \rho_t^{r\varepsilon} \sup_{|x|\leq r} \left| \frac{\xi_t S_\varepsilon \theta_x f}{\xi_t B_{\varepsilon x}^\varepsilon} - \lambda^2 f \right|$$

$$\leq E_\mu \rho_t^{r\varepsilon} \left| \frac{\xi_t S_\varepsilon f_r^+}{\xi_t B_0^{\varepsilon(1-r)}} - \frac{\lambda^2 f_r^+}{(1-r)^2} \right| + E_\mu \rho_t^{r\varepsilon} \left| \frac{\xi_t S_\varepsilon f_r^-}{\xi_t B_0^{\varepsilon(1+r)}} - \frac{\lambda^2 f_r^-}{(1+r)^2} \right|$$

$$+ E_\mu \rho_t^{r\varepsilon} \lambda^2 ((1-r)^{-2} f_r^+ - (1+r)^{-2} f_r^-).$$

Dividing by $E_\mu \rho_t^{r\varepsilon}$ and applying (24) to $f_r^\pm$, we get as $\varepsilon \to 0$

$$\limsup_{\varepsilon\to 0} E_\mu^0(|\xi_t S_\varepsilon f/\rho_t^\varepsilon - \lambda^2 f| \wedge 1)$$
$$\leq \sup_{|x|\leq r\varepsilon} \|P_\mu^0 - P_{\mu\theta_x}^0\|_C + \lambda^2((1-r)^{-2} f_r^+ - (1+r)^{-2} f_r^-),$$

for any neighborhood $C$ of 0. Here both terms on the right tend to 0 as $r \to 0$, the former by Lemma 3.5 and the latter by the continuity of $f$ and dominated convergence. Hence,

$$\xi_t S_\varepsilon f/\rho_t^\varepsilon \xrightarrow{d} \lambda^2 f \text{ under } P_\mu^0, \qquad f \in C_K^2,$$

and (i) follows by Theorem 16.16 in [13].

(iii) Letting $\varepsilon \to 0$ for fixed $x \in \mathbb{R}^d$ and $r > 0$, we get by Theorem 5.3(ii) and Lemmas 3.2(iv) and 5.1

$$\frac{P_\mu\{\xi_t B_{\varepsilon x}^{\varepsilon r} > 0\}}{P_\mu\{\xi_t B_0^\varepsilon > 0\}} \sim \frac{m(r\varepsilon)}{m(\varepsilon)} \frac{|\log \varepsilon|}{|\log(r\varepsilon)|} \frac{(\mu * p_t)(\varepsilon x)}{\mu p_t} \to 1.$$

Keeping $x$ and $r$ fixed and choosing $c > 0$ with $B_x^r \subset B_0^c$, we get in particular $P_\mu[\xi_t B_0^\varepsilon > 0|\xi_t B_0^{c\varepsilon} > 0] \to 1$, and so as $\varepsilon \to 0$

$$P_\mu[\xi_t B_{\varepsilon x}^{\varepsilon r} > 0|\xi_t B_0^\varepsilon > 0]$$
$$\geq \frac{P_\mu\{\xi_t B_{\varepsilon x}^{\varepsilon r} > 0\}}{P_\mu\{\xi_t B_0^\varepsilon > 0\}} - \frac{P_\mu(\{\xi_t B_0^\varepsilon > 0\}\Delta\{\xi_t B_0^{c\varepsilon} > 0\})}{P_\mu\{\xi_t B_0^\varepsilon > 0\}} \to 1.$$

The assertion follows since $x$ and $r$ were arbitrary.

(ii) For any $f \in C_K^2$, we get by (23), (iii), Lemma 4.2(ii) and Jensen's inequality

$$(E_\mu^\varepsilon |\xi_t S_\varepsilon f - \rho_t^\varepsilon \lambda^2 f|)^2 \leq E_\mu^\varepsilon (\xi_t S_\varepsilon f - \rho_t^\varepsilon \lambda^2 f)^2$$
$$\lesssim \frac{E_\mu (\xi_t S_\varepsilon f - \rho_t^\varepsilon \lambda^2 f)^2}{P_\mu \{\rho_t^\varepsilon > 0\}} \ll \varepsilon^4 |\log \varepsilon|^2.$$

Similarly, we see from Lemmas 3.2 and 4.2(ii) that

$$E_\mu \rho_t^\varepsilon = \frac{E_\mu \rho_t^\varepsilon}{P_\mu \{\rho_t^\varepsilon > 0\}} \asymp \varepsilon^2 |\log \varepsilon|.$$

The result follows by combination of these estimates. □

We may finally use the results of Section 8 to show that the local invariance fails for $d \geq 3$. The argument also shows that the main results of Section 8 have no counterparts for $d = 2$.

PROPOSITION 9.2. *For $d \geq 3$, let the DW-process $\xi$ in $\mathbb{R}^d$ be locally finite under $P_\mu$, and fix any $t > 0$. Letting $\varepsilon \to 0$ and then $h \to 0$, we have:*

(i) $P_\mu[\xi_t S_\varepsilon B_x^h = 0 | \xi_t B_0^\varepsilon > 0] \to 1$ *for all $x \in \mathbb{R}^d$,*
(ii) $P_\mu^0 \{\xi_t S_\varepsilon B_x^h = 0\} \to 1$ *for all $x \neq 0$ in $\mathbb{R}^d$.*

PROOF. For any bounded initial measure $\mu$ on $\mathbb{R}^d$, we have $\lambda^d(\operatorname{supp} \xi_t) = 0$ a.s., for example, by Theorem 7.1(i). Using Fubini's theorem (to ensure measurability) and Theorem 8.1(i), we get for any $B \in \hat{\mathcal{B}}^d$

$$0 = \varepsilon^{2-d} P_\mu \{\lambda^d (\varepsilon B \cap \operatorname{supp} \xi_t) > 0\}$$
$$\to \mu p_t \tilde{P} \{\lambda^d (B \cap \operatorname{supp} \tilde{\xi}) > 0\},$$

which implies $\lambda^d(\operatorname{supp} \tilde{\xi}) = 0$ a.e. $\tilde{P}$. By the stationarity of $\tilde{\xi}$ and the shift invariance of the function $\lambda^d(\operatorname{supp} \mu)$, the same property holds a.s. under $\tilde{P}^0$.

Next, Fubini's theorem yields $\tilde{P}^0 \{x \in \operatorname{supp} \tilde{\xi}\} = 0$ for $x \in \mathbb{R}^d$ a.e. $\lambda^d$. In particular, we may choose an $x \neq 0$ with $x \notin \operatorname{supp} \tilde{\xi}$ a.e. $\tilde{P}^0$. By rotational symmetry and scaling invariance, this remains true for every $x \neq 0$. Since $\operatorname{supp} \tilde{\xi}$ is closed, Theorem 8.2 yields

$$\lim_{h \to 0} \limsup_{\varepsilon \to 0} P_\mu^0 \{\xi_t S_\varepsilon B_x^h > 0\} = \lim_{h \to 0} \tilde{P}^0 \{\tilde{\xi} B_x^h > 0\} = 0,$$

proving (ii). Assertion (i) holds by a similar argument based on Theorem 8.1(i). □

**Acknowledgment.** I am grateful to a referee for pointing out some interesting connections to results in [20, 23, 24]. I am also grateful to Jean-François Le Gall for pointing out that part (i) of Theorem 7.1 was obtained by Tribe [27] (see also [6, 25, 27]). However, the result for $d = 2$ seems to be new. Our proof uses a new approach, part of which is explained in detail only for $d \geq 3$. Thus, the entire proof of Theorem 7.1 is still needed.

DEPARTMENT OF MATHEMATICS AND STATISTICS
AUBURN UNIVERSITY
221 PARKER HALL
AUBURN, ALABAMA 36849
USA
E-MAIL: kalleoh@auburn.edu